\documentclass[oneside,english]{siamltex}
\usepackage[T1]{fontenc}
\usepackage[latin9]{inputenc}
\usepackage{geometry}
\usepackage{array}
\usepackage{float}
\usepackage{dvipost}
\usepackage{multirow}
\usepackage{amsbsy}
\usepackage{amssymb}
\usepackage{graphicx}
\usepackage{esint}
\usepackage[numbers]{natbib}

\makeatletter

\providecommand{\tabularnewline}{\\}
\dvipostlayout
\dvipost{osstart color push Red}
\dvipost{osend color pop}
\dvipost{cbstart color push Blue}
\dvipost{cbend color pop}

\newcommand\eqref[1]{(\ref{#1})}

\usepackage{lineno}
\usepackage{tikz}
\usetikzlibrary{arrows,decorations.pathreplacing}
\newlength{\ypos}
\newtheorem{conjecture}{Conjecture}[section]
\newtheorem{rmk}{Remark}[section]

\@ifundefined{showcaptionsetup}{}{%
 \PassOptionsToPackage{caption=false}{subfig}}
\usepackage{subfig}
\makeatother

\usepackage{babel}

\begin{document}

\title{Numerical simulation of solutions and moments of the smoluchowski
coagulation equation}

\author{D.~D.~Keck\thanks{Department of Applied Mathematics, University of Colorado, Boulder, CO 80309-0526}
\and D.~M.~Bortz$^{*}$\thanks{BioFrontiers Institute, University of Colorado, Boulder, CO 80309-0596}}
\maketitle
\begin{abstract}
Researchers have employed variations of the Smoluchowski coagulation
equation to model a wide variety of both organic and inorganic phenomena
and with relatively few known analytical solutions, numerical solutions
play an important role in studying this equation. In this article,
we consider numerical approximations, focusing on how different discretization
schemes impact the accuracy of approximate solution moments. Pursuing
the eventual goal of comparing simulated solutions to experimental
data, we must carefully choose the numerical method most appropriate
to the type of data we attain. Within this context, we compare and
contrast the accuracy and computational cost of a finite element approach
and a finite volume-based scheme.

Our study provides theoretical and numerical evidence that the finite
element approach achieves much more accuracy when the system aggregates
slowly, and it does so with much less computation cost. Conversely,
the finite volume method is slightly more accurate approximating the
zeroth moment when the system aggregates quickly and is much more
accurate approximating the first moment in general. 

Lastly, our study also provides numerical evidence that the finite
element method (conventionally considered first order) actually belongs
to a class of discontinuous Galerkin methods that exhibit superconvergence,
or second order in our case.
\end{abstract}

\section{Introduction}

The Smoluchowski coagulation equation was originally developed by
Marian von Smoluchowski in the early 1900's \citep{Smoluchowski1916,Smoluchowski1917}
to study gelling colloids. More recently, researchers have employed
variations of this model to study organic phenomena such as bacterial
growth \citep{BortzEtal2008bmb}, marine snow \citep{Kiorboe2001},
algal blooms \citep{Ackleh1997NATMA,AcklehFitzpatrick1997,Riebesell1992},
and schooling fish \citep{Niwa1998} and inorganic phenomena such
as powder metallurgy \citep{Kumar2009}, astronomy \citep{Lee2000icarus,Lee2001jopa,Makino1998,Silk1978},
aerosols \citep{Drake1972}, irradiation of metals \citep{Surh2008},
and meteorology \citep{Pruppacher1980}. Furthermore, because only
a few known analytic solutions exist, numerical solutions to the Smoluchowski
coagulation equation play an important role in the study of this equation
\citep{Lee2001jopa}. A number of computational approaches have been
formulated using finite elements \citep{BortzEtal2008bmb,Mahoney2002},
finite volumes \citep{Filbet2004,Verkoeijen2002}, successive approximations
\citep{Ramkrishna2000}, method of moments \citep{Barrett1996,Madras2004},
Monte Carlo simulations \citep{Kruis2000,Lin2002}, and mesh-free
approaches that capitalize on radial basis functions \citep{Ranjbar2010}.
This overview is just a survey of some important works in the field.
For more comprehensive treatments we direct the interested reader
to the reviews concerning mathematical results by Wattis \citep{Wattis2006a}
and Aldous \citep{Aldous1999} and a an empirical comparison of numerical
techniques by Lee \citep{Lee2001jopa}.

Before proceeding, we note that while a mean-field approach to modeling
particles in suspension can be very useful, we must use care with
our terminology. Throughout this paper, we refer to two types of \emph{distributions}
over a particle volume domain of interest, ${\bf X}\subset\mathbb{R}_{+}$.
In both cases, we use distribution in the sense that we identify a
quantity of aggregates per total volume of the aggregates in ${\bf X}$.
First we denote a \emph{size distribution} $f=f(t,x)\ge0$ as the
\textit{number} density of aggregates of a given volume $x$ at time
$t\ge0$. We denote a \emph{volume distribution} $g=g(t,x)$ as the
\textit{volume} density of aggregates of a given volume $x$. A superscript
$N$\emph{ }will denote numerical approximation, e.g., $f^{N}$, \emph{$g^{N}$,
etc}. Note that the volume distribution relates to the size distribution
as $g(t,x)=xf(t,x)$.

To further promote clarity, we define the \textit{partial $i^{th}$
moment} as
\begin{equation}
M_{i}(f(t,\cdot);x_{1},x_{2})=\int_{x_{1}}^{x_{2}}x^{i}f(t,x)dx\label{eq:ith moment}
\end{equation}
where $f(t,\cdot)$ is the size distribution at time $t$. For example,
$M_{0}(f(t,\cdot);x_{1},x_{2})$ or $M_{0}(\frac{g(t,\cdot)}{(\cdot)};x_{1},x_{2})$
represents the total number of clusters having volumes between $x_{1}$
and $x_{2}$, and $M_{1}(f(t,\cdot);x_{1},x_{2})$ or $M_{1}(\frac{g(t,\cdot)}{(\cdot)};x_{1},x_{2})$
represents the total volume of clusters where each cluster included
in the total volume individually has a volume between $x_{1}$ and
$x_{2}$.

When using models derived from the Smoluchowski equation to study
the real world, we will compare a simulated solution to experimental
data with the eventual goal to illuminate some scientific phenomenon.
Naturally, different experiments yield different types of data, just
as different discretization schemes yield varying accuracies for the
approximations of different quantities. Accordingly, the type of data
should guide the choice of numerical scheme. For example, experimentalists
utilizing a Coulter counter \citep{Zhe2007} often provide data in
the form of a vector of partial zeroth moments. Conversely, when employing
a flow cytometer \citep{Shapiro2005} or dynamic light scattering
instrumentation \citep{Berne2000}, data is often reported as a partial
first moment. Therefore it is important to choose a discretization
scheme based on which moment is reported by the specific experimental
apparatus. In addition, no experimental device will provide full information
about the whole positive real axis; there will always be limited ranges
for reliable data. Accordingly, we must also address the additional
issue of how to deal with a lack of information about particle aggregates
outside the detection limits of a given device.

In this paper, we consider a \textit{Finite Element Method} (\textit{FEM})
approach developed in Banks and Kappel \citep{BanksKappel1989} (extended
by Ackleh and Fitzpatrick \citep{AcklehFitzpatrick1997}, and explored
in Bortz et al.~\citep{BortzEtal2008bmb}). We also consider a finite
volume-type scheme, which we designate as the \textit{Filbet and Laurençot
Flux Method (FLFM),} developed in \citep{Filbet2004}. 

For both discretization approaches, we pay particular attention to
the aggregate volume domain $\mathbf{X}$ and its limits. With both
discretizations of ${\bf X},$ we lose information, and the impact
of that loss on the respective method's accuracies deserves investigation.
In that light, a goal of this work is to fully compare the two schemes
in terms of their accuracy in approximating a solution and their accuracy
in approximating zeroth and first moments and in terms of computation
cost. Our investigation supports second order convergence of both
methods to fine grid solutions and to the zeroth and first moments.
Additionally, our investigation reveals that when modeling slowly
aggregating systems, the FEM can provide as little error approximating
a true solution as the FLFM and more accurately approximates the zeroth
moment, but does so with significantly less computation cost. Conversely,
the FLFM approximates the zeroth moment slightly more accurately for
slowly aggregating systems, and approximates the first moment more
accurately in general.

The paper is organized as follows. Section \ref{sec:Numerical-Methods}
describes each of the two methods. Section \ref{sec:Analytical-Solutions}
provides a description of the explicit solutions studied. Section
\ref{sec:Computational-Results} discusses the results of our numerical
computations. Finally, Section \ref{sec:Conclusions-and-Future} summarizes
the conclusions of this work.

\section{Numerical Methods\label{sec:Numerical-Methods}}

Of the multiple numerical schemes mentioned in the introduction, we
restrict ourselves to the FEM originally developed in Banks and Kappel
\citep{BanksKappel1989} and extended by Ackleh and Fitzpatrick \citep{AcklehFitzpatrick1997}
and the FLFM described by Filbet and Laurençot in \citep{Filbet2004}.
In Section \ref{sub:Num Meth Gen Descrip}, we give a brief description
of the Smoluchowski coagulation equation and the common assumptions
both methods use in numerically solving the equation. In Section \ref{sub:Num Meth FEM},
we highlight the most important parts of the FEM model and the aggregation
vectors it creates. Then in Section \ref{sub:Num Meth FLFM}, we do
the same for the FLFM model.

\subsection{\label{sub:Num Meth Gen Descrip}Model and Discretization Overview}

In the early 1900's, van Smoluchowski developed a model to study the
coagulation of colloids,

\begin{equation}
\frac{d}{dt}f_{k}=\frac{1}{2}\sum_{i+j=k}K(i,j)f_{i}f_{j}-\sum_{i}K(i,k)f_{i}f_{k},
\end{equation}
where $f_{k}$ represents the number density of aggregates of volume
$k$, and $K(i,j)$ is the aggregation kernel denoting the rate at
which aggregates of size $i$ and $j$ form a combined aggregate of
size $i+j$ \citep{BortzEtal2008bmb,Smoluchowski1916,Smoluchowski1917}.
Müller subsequently extended this model to a continuous PDE \citep{Muller1928,Filbet2004} 

\begin{eqnarray}
\partial_{t}f & = & A(f),\,\,(t,x)\in\mathbb{R}_{+}^{2},\label{eq:smolu}\\
f(0) & = & f_{0},\,\, x\in\mathbb{R}_{+}\nonumber 
\end{eqnarray}
 where we describe each aggregate solely by its volume $x>0$, with
$f=f(t,\cdot)$ representing the continuous size distribution function
of aggregates at time $t\ge0$. The coagulation term is

\begin{eqnarray}
A(f) & = & A_{in}(f)-A_{out}(f)\nonumber \\
 & = & \frac{1}{2}\int_{0}^{x}K_{A}(y,x-y)f(t,y)f(t,x-y)dy-f(x)\int_{0}^{\infty}K_{A}(x,y)f(y)dy\label{eq:contin agg}
\end{eqnarray}
where $K_{A}(x,y)$ is the aggregation kernel indicating the rate
at which aggregates of volumes $x$ and $y$ join together creating
an aggregate of volume $x+y$. Notice the first integral, $A_{in}(f)$,
describes aggregates with volumes $y$ and $x-y$ aggregating to a
combined volume $x$, and the second integral, $A_{out}(f)$, models
interactions between the aggregate of volume $x$ with all other aggregates
of volume $y$ forming an aggregate of volume $x+y$. Also, note that
the aggregation kernel $K_{A}(x,y)$ is positive and symmetric
\[
0<K_{A}(x,y)=K_{A}(y,x),\,(x,y)\in\mathbb{R}_{+}^{2},
\]
as well as homogeneous, which literature in this field defines as
\begin{equation}
K_{A}(\lambda x,\lambda y)=\lambda^{m}K_{A}(x,y),\,\,\,\lambda>0,\, m\ge0,\, x,y<\infty.\label{eq:homgen def}
\end{equation}
Because only aggregation is considered, the total number of particles
decreases with each coagulation event. 

The conservation properties of the model in (\ref{eq:contin agg})
warrant a brief discussion. The model is based on conservation of
mass principles, but practically speaking, in our simulations, we
expect to \textit{lose} mass from the system. First, for aggregation
kernels such as the multiplicative kernel, $K_{A}(x,y)=xy$, the system
experiences growth rapid enough that aggregates with \textit{infinite}
volume develop in finite time \citep{Wattis2006a}. This phenomenon
is commonly referred to as \textit{gelation}. Mass is not physically
lost, but the aggregates with \textit{infinite} volume possess fundamentally
different mathematical properties than the individual aggregates that
make up the gel. We direct the interested reader to \citep{Ziff1980},
in which Ziff and Stell provide a thorough description of the implications
of various assumptions on the post-gelation behavior of the solutions
and of the moments. Second, when we discretize the equation, the domain
will naturally have a finite maximum value, $x_{max}$, replacing
the infinite upper integration limit of the second integral in (\ref{eq:contin agg}).
Therefore, the approximate solution will not include any impact of
aggregates larger than $x_{max}$. 

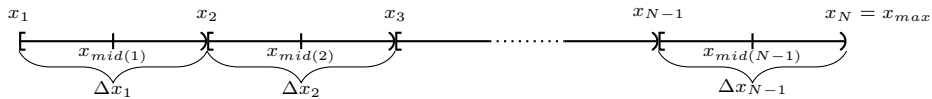
\begin{figure}[H]
\centering{}\setlength{\ypos}{-0.8cm} 
\begin{tikzpicture}[scale=1]     
	\begin{scope}[thick,font=\scriptsize]
		\draw [-] (0,0) -- (0,0) node [above=4pt] {$x_1$}; 	
		\draw [-] (1.25,-3pt) -- (1.25,3pt) node [below=2pt] {$x_{mid(1)}$}; 	
		\node (x1) at (0,-3pt) {}; 	
		\node (x2) at (2.5,-3pt) {}; 		\draw[decorate,decoration={brace,mirror,amplitude=12pt,raise=-3pt},thin] 		(x1.south)-- node [below=5pt] {$\Delta x_1$} (x2.south);
		\draw [[-)] (0,0) -- (2.5,0) node [above=4pt] {$x_2$}; 	
		\draw [-] (3.75,-3pt) -- (3.75,3pt) node [below=2pt] {$x_{mid(2)}$}; 	\draw [[-)] (2.5,0) -- (5,0) node [above=4pt] {$x_3$}; 	
		\node (x3) at (5,-3pt) {}; 		\draw[decorate,decoration={brace,mirror,amplitude=12pt,raise=-3pt},thin] 		(x2.south)-- node [below=5pt] {$\Delta x_2$} (x3.south);
		\draw [[-] (5,0) -- (6.25,0); 	\draw [dotted] (6.25,0) -- (7.25,0); 	
		\draw [-)] (7.25,0) -- (8.5,0) node [above=4pt] {$x_{N-1}$}; 	
		\draw [[-)] (8.5,0) -- (11,0) node [above=4pt] {$~~~~~~~~~x_N=x_{max}$}; 	
		\draw [-] (9.75,-3pt) -- (9.75,3pt) node [below=2pt] {$x_{mid(N-1)}$}; 	
		\node (xNm1) at (8.5,-3pt) {}; 	\node (xN) at (11,-3pt) {}; 		\draw[decorate,decoration={brace,mirror,amplitude=12pt,raise=-3pt},thin] 		(xNm1.south)-- node [below=4pt] {$\Delta x_{N-1}$} (xN.south);
	\end{scope} 
\end{tikzpicture}\caption{Discretization - Grid used for both FEM and FLFM where $\Delta x=\Delta x_{i}=x_{i+1}-x_{i\,}\forall i.$
Only uniform grids are used in this study.}
\label{Discret}
\end{figure}

In the study that follows, we consider two aggregation kernels, the
identity kernel, $K_{A}(x,y)\equiv1$, and the multiplicative kernel,
$K_{A}(x,y)=xy$ (both of which have experienced widespread use).
In order to solve the governing PDE, (\ref{eq:smolu}), we must first
discretize ${\bf X}$, which leads to a system of ODEs. We then advance
the solution to the system of ODEs in discrete time steps by employing
a variable order solver based on the numerical differentiation formulas
(NDFs) as implemented in Matlab's \textbf{ode15s}. In this study,
we only consider uniform grids, so $\Delta x=\Delta x_{i}=x_{i+1}-x_{i,\,}\forall i$.
To provide a sense for the subtle differences between the grids used
with both approaches, as well as to provide a visual representation
of the notation used throughout the rest of this study, we have included
Figure \ref{Discret}. Note that $x_{mid(i)}$ is the midpoint between
$x_{i}$ and $x_{i+1}$. 

\subsection{\label{sub:Num Meth FEM}Finite Element Approach (FEM) }

\begin{flushleft}
We now provide a brief overview of the important details of the FEM
as discussed in Bortz et al.~\citep{BortzEtal2008bmb} and Ackleh
and Fitzpatrick \citep{AcklehFitzpatrick1997}. First we define our
solution space as $H=L^{1}([x_{min},x_{max}],\mathbb{R}_{+})$, the
space of integrable functions mapping a closed, bounded subset of
positive reals into the positive reals where $x_{min}$ and $x_{max}$
are the respective minimum and maximum aggregate volumes. We use $N-1$
elements in our numeric grid (see Figure \ref{Discret}), where each
element boundary is denoted $x_{i}$ for $i=1,\dots,N$ such that
$x_{1}=x_{min}$ and $x_{N}=x_{max}$. We then choose a set of hat
basis functions for $i\in[1,N-1]$ 
\[
\beta_{i}^{N}(x)=\left\{ \begin{array}{ll}
1; & x_{i}^{N}\le x<x_{i+1}^{N}\\
0; & \mbox{otherwise}
\end{array}\right\} 
\]
that form an orthogonal basis for our solution space
\[
H^{N}=\left\{ h\in H:\, h=\sum_{i=1}^{N-1}\alpha_{i}\beta_{i}^{N},\,\alpha_{i}\in\mathbb{R}\right\} .
\]
The coefficients are $\alpha_{j}=\frac{1}{\Delta x}\int_{x_{j}}^{x_{j+1}}h(x)dx$,
which allows us to define the projections $\pi^{N}:H\rightarrow H^{N}$
as
\[
\pi^{N}h=\sum_{j=1}^{N-1}\alpha_{j}\beta_{j}^{N}
\]
providing an orthogonal projection of $H$ onto $H^{N}$. We can now
approximate (\ref{eq:smolu}) with a set of $N-1$ ODEs 
\begin{eqnarray}
f_{t}^{N} & = & \pi^{N}\left(A(f^{N})\right),\label{eq:FEMdiscrete}\\
f^{N}(0,x) & = & \pi^{N}\left(f(0,x)\right).\nonumber 
\end{eqnarray}
The discretization described in (\ref{eq:FEMdiscrete}) is a discontinuous
Galerkin (DG) approximation using a zeroth order polynomial basis,
and from Section 4 of \citep{AcklehFitzpatrick1997}, $f^{N}\rightarrow f$
uniformly in norm. As reported in \citep{AcklehFitzpatrick1997},
the FEM converges at first order on $L^{\infty}[\mathbf{X}]$. However,
recent work by Yang and Shu in \citep{YangShu2012sinum} on the supercovergence
of DG methods and results from our numerical experiments support second
order convergence in $L^{1}[\boldsymbol{X}]$. Future efforts will
include an analytical investigation of the following conjecture.
\par\end{flushleft}
\begin{conjecture}
The approximate solution, $f^{N}(t_{k},x)$, converges to the analytical
solution with order 2 in $L^{1}[\boldsymbol{X}]$.
\end{conjecture}
\begin{flushleft}
Using this result, we note that the convergence of the $i$th moment
is second order (or faster) in $\Delta x$.
\par\end{flushleft}

\begin{rmk}
\it{The i$^{th}$ moment $M_i^N$ of the approximate solution converges to the i$^{th}$ moment $M_i$ of the analytical solution with order 2 (or faster) in $L^1[\bold{X}]$.}
\end{rmk}
\begin{proof}
Recall that for fixed $t_{k}$, there is second order convergence
in $\Delta x$ of $f^{N}$ to $f$. 
\begin{eqnarray*}
\lim_{N\rightarrow\infty}\left|M_{i}^{N}(f(t_{k},\cdot);\mathbf{X})-M_{i}(f(t_{k},\cdot);\mathbf{X})\right| & = & \lim_{N\rightarrow\infty}\left|\int_{\mathbf{X}}x^{i}f^{N}(t_{k},x)dx-\int_{\mathbf{X}}x^{i}f(t_{k},x)dx\right|\\
 & = & \lim_{N\rightarrow\infty}\left|\int_{\mathbf{X}}x^{i}\left\{ f^{N}(t_{k},x)-f(t_{k},x)\right\} dx\right|\\
 & \le & \lim_{N\rightarrow\infty}\int_{\mathbf{X}}x^{i}\left|f^{N}(t_{k},x)-f(t_{k},x)\right|dx\\
 & \le & \int_{\mathbf{X}}x^{i}dx\lim_{N\rightarrow\infty}\left\Vert f^{N}(t_{k},x)-f(t_{k},x)\right\Vert _{L^{1}[x_{j},x_{j+1}]}
\end{eqnarray*}

\end{proof}
\begin{flushleft}
Before giving a detailed description of each term in the discretized
system, we must highlight a few considerations. For the FEM, we use
$f_{i}^{N}$ for $i\in[1,N-1]$ to represent the discretized size
distribution of aggregates in the $i^{th}$ element. In our discretization
scheme, the $i^{th}$ element spans $[x_{i},x_{i+1})$ with $x_{1}=x_{min}$
(depicted in Figure \ref{Discret}). For simplification of the explanations
below, we let $x_{1}=0$. In the validation section, however, we compare
our solution with an analytical solution to (\ref{eq:smolu}), which
has a singularity at zero. Therefore in practice, we set $x_{1}>0$
in some scenarios. We then take the approach that all aggregates in
a given element, $i$, have the discrete volume infinitesimally close
to $x_{i+1}$ with discretized size distribution, $f_{i}^{N}$. Under
that rationale, no two aggregates can combine to produce one of size
$x_{2}$, so the rate of change of $f_{1}^{N}$ is strictly negative.
Furthermore, we account for the interactions of all particles with
volumes up to $x_{max}$. As mentioned above, the implication of this
choice is that when two particles with volumes $x_{i}$ and $x_{j}$
aggregate to form a particle with volume $x_{i}+x_{j}$ and $x_{i}+x_{j}>x_{N}=x_{max}$,
we lose mass from our system. 
\par\end{flushleft}

\begin{flushleft}
We can now describe our discretized system fully. At $t=0$, we project
our initial conditions such that,
\begin{equation}
\pi^{N}(f(0,\cdot))=\sum_{j=1}^{N-1}\alpha_{j}\beta_{j}^{N}(\cdot)\label{eq:init cond FEM}
\end{equation}
and for $x\in[x_{i},x_{i+1})$
\[
\sum_{j=1}^{N-1}\alpha_{j}\beta_{j}^{N}(x)=\frac{1}{\Delta x}\int_{x_{i}}^{x_{i+1}}f(0,y)dy.
\]
Then for each discrete time step, $t_{k}$, we create two discretized
vectors with $N-1$ elements. The first vector represents aggregation
of particles out of each element, while the second represents aggregation
of particles into each element. These two vectors take on a different
form for each of the two aggregation kernels ($K_{A}(x,y)\equiv1$
and $K_{A}(x,y)=xy$) used in our study, and we present both vectors
below. For aggregates which coagulate with an aggregate of size $x$
to aggregates greater than size $x$, i.e., \textit{aggregation out},
\[
A_{out}(f)=-f(x)\int_{0}^{\infty}K_{A}(x,y)f(y)dy.
\]
To discretize this integral at each $t_{k}$, we have to truncate
it to some finite maximum, $x_{max}$. For the general aggregation
kernel, $K_{A}(x,y)$, and for a given element, $i$,
\begin{eqnarray}
\pi^{N}\left(A_{out}(f_{i}^{N})\right) & = & -\alpha_{i}\int_{x_{1}}^{x_{max}}K_{A}(x_{i+1},y)\alpha(y)dy\\
 & = & -\alpha_{i}\sum_{j=1}^{N-1}\left(\int_{x_{j}}^{x_{j+1}}K_{A}(x_{i+1},y)\alpha(y)dy\right).\label{eq:FEMaggoutintegrals}
\end{eqnarray}
Now note that $\alpha_{j}=\frac{1}{\Delta x}\int_{x_{j}}^{x_{j+1}}h(x)dx=\frac{1}{\Delta x}\int_{x_{j}}^{x_{j+1}}f_{j}^{N}dx=f_{j}^{N}$
for any $t_{k}$ where $k>0$, therefore $\alpha_{j}$ is simply the
approximated, discrete size distribution, $f_{j}^{N}$, from the previous
time step. Then after making the appropriate substitutions and integrating,
for $K_{A}(x,y)\equiv1$,
\begin{equation}
\pi^{N}\left(A_{out}(f_{i}^{N})\right)=-f_{i}^{N}\Delta x\sum_{j=1}^{N-1}f_{j}^{N},\label{eq:ais1out}
\end{equation}
and for $K_{A}(x,y)=xy$,
\begin{equation}
\pi^{N}\left(A_{out}(f_{i}^{N})\right)=\frac{-f_{i}^{N}x_{i+1}}{2}\left[\sum_{j=1}^{N-1}\left(x_{j+1}^{2}-x_{j}^{2}\right)f_{j}^{N}\right].\label{eq:aisxyout}
\end{equation}
 We enter (\ref{eq:ais1out}) and (\ref{eq:aisxyout}) for our respective
\textit{aggregation out} vectors. Now we consider the aggregates which
coagulate to form an aggregate of size $x$, i.e., \textit{aggregation
in}. In the continuous case,
\[
A_{in}(f)=\frac{1}{2}\int_{0}^{x}K_{A}(y,x-y)f(t,y)f(t,x-y)dy.
\]
Discretizing this integral for the general aggregation kernel, $K_{A}(x,y)$,
and for a given element, $i>1$,
\begin{eqnarray}
\pi^{N}\left(A_{in}(f_{i}^{N})\right) & = & \frac{1}{2}\int_{x_{1}}^{x_{i}}K_{A}(y,x_{i}-y)\alpha(y)\alpha(x_{i}-y)dy\label{eq:FEMaisxy_1}\\
 & = & \frac{1}{2}\sum_{j=1}^{i-1}\int_{x_{j}}^{x_{j+1}}K_{A}(y,x_{i}-y)\alpha_{j}\alpha(x_{i}-y)dy\label{eq:FEMaisxy_2}
\end{eqnarray}
Again note that $\alpha_{j}=\frac{1}{\Delta x}\int_{x_{j}}^{x_{j+1}}h(x)dx=\frac{1}{\Delta x}\int_{x_{j}}^{x_{j+1}}f_{j}^{N}dx=f_{j}^{N}$
for any $t_{k}$ where $k>0$, therefore $\alpha_{j}$ is simply the
approximated, discrete size distribution, $f_{j}^{N}$, from the previous
time step. Then after making the appropriate substitutions and integrating
for $K_{A}(x,y)\equiv1$,
\begin{equation}
\pi^{N}\left(A_{in}(f_{i}^{N})\right)=\frac{1}{2}\Delta x\sum_{j=1}^{i-1}f_{j}^{N}f_{i-j}^{N},\label{eq:kais1in}
\end{equation}
and for $K_{A}(x,y)=xy$,
\begin{equation}
\pi^{N}\left(A_{in}(f_{i}^{N})\right)=\frac{1}{2}\sum_{j=1}^{i-1}\left[x_{i}\frac{x_{j+1}^{2}-x_{j}^{2}}{2}-\frac{x_{j+1}^{3}-x_{j}^{3}}{3}\right]f_{j}^{N}f_{i-j}^{N}.\label{eq:aisxyin}
\end{equation}
We enter (\ref{eq:kais1in}) and (\ref{eq:aisxyin}) for our respective
\textit{aggregation in} vectors. 
\par\end{flushleft}

\subsection{\label{sub:Num Meth FLFM}Filbet and Laurençot Flux Method (FLFM)}

Our second approach for numerically solving the Smoluchowski coagulation
equation is a scheme developed by Filbet and Laurençot in \citep{Filbet2004}.
They base their scheme on a finite volume method and calculate a mass
flux quantity, which we denote, $J[f](t,x)$. Here, $J[f](t,x)$ represents
a mass flux from aggregates with volumes at most $x$ to aggregates
with volumes greater than $x$. Then similar to Filbet and Laurençot,
we reformulate \citep{Filbet2004,Makino1998} the Smoluchowski coagulation
equation (\ref{eq:smolu})
\begin{equation}
x\partial_{t}f=-\partial_{x}J[f],\,(t,x)\in\mathbb{R}_{+}^{2}\label{eq:smolu_flux}
\end{equation}
where
\begin{equation}
J[f](t,x)=\int_{0}^{x}\int_{x-u}^{\infty}uK_{A}(u,v)f(t,u)f(t,v)dvdu,\, x\in\mathbb{R}_{+};\, t\in\mathbb{R}_{+}.\label{eq:contin flux}
\end{equation}
We can recover the conventional formulation of the Smoluchowski coagulation
equation, (\ref{eq:smolu}), by substituting (\ref{eq:contin flux})
into (\ref{eq:smolu_flux}) and applying Leibniz's rule. Once again,
for numerical purposes, we must truncate the volume variable, $x$,
to a finite value, $x_{max}$. Filbet and Laurençot \citep{Filbet2004}
describe a number of choices (conservative and nonconservative) for
truncating the inner integral in (\ref{eq:contin flux}). We follow
their recommendation by using the nonconservative formulation
\begin{equation}
J[f](t,x)=\int_{0}^{x}\int_{x-u}^{x_{max}}uK_{A}(u,v)f(t,u)f(t,v)dvdu,\, x\in(0,x_{max});\, t\in\mathbb{R}_{+}\label{eq:Jnoncon-1}
\end{equation}
because it generates accurate approximations when solutions include
gelation \citep{Bak1991,Costa1998,Filbet2004}. As mentioned in \citep{Filbet2004},
we can rewrite (\ref{eq:smolu_flux}) in terms of a volume distribution,
which we denote $g(t,x)=xf(t,x)$. Then (\ref{eq:smolu_flux}) and
(\ref{eq:contin flux}) become
\begin{equation}
\partial_{t}g=-\partial_{x}J[g],\,(t,x)\in\mathbb{R}_{+}^{2}\label{eq:smolu_flux_g}
\end{equation}
and
\begin{equation}
J[g](t,x)=\int_{0}^{x}\int_{x-u}^{x_{max}}\frac{K_{A}(u,v)}{v}g(t,u)g(t,v)dvdu,\, x\in(0,x_{max});\, t\in\mathbb{R}_{+}.\label{eq:contin flux_g}
\end{equation}
This formulation is especially useful in application when the data
has the form of a volume distribution, and we discuss the advantages
of this formulation in more detail in Section \ref{sub:Moment calcs}.

Numerically solving the formulation in (\ref{eq:smolu_flux}) differs
fundamentally from numerically solving the formulation for which we
used finite elements. With FEM, we track the changes in size distribution
for given elements. With the formulation in (\ref{eq:smolu_flux}),
we track the discretized volume distribution, which we denote, $g_{i}^{N}(t_{k})$,
representing the approximated mean value of $g(t_{k},x)$ in the element,
$[x_{i},x_{i+1})$ at discrete time steps, $t_{k}$. This formulation
also includes the discretized mass flux, which we denote, $J_{i}^{N}(t_{k})$,
across element boundaries (recall Figure \ref{Discret} with boundaries,
$x_{i}$) at discrete time steps. To understand how we determine the
flux, $J_{i}^{N}(t_{k}),$ at each element boundary, $x_{i}$, for
any time step, consider the following. First, $J_{1}^{N}(t_{k})=0\,\forall k$,
and for all other element boundaries, flux across a given boundary,
$x_{r}$, requires $x_{mid(i)}+\tilde{x}\ge x_{r}$. We now fix $x_{r}$
and $x_{mid(i)}<x_{r}$ to determine the discretized volumes, $x_{mid(j)}$,
that are equal to or larger than $\tilde{x}$. For each $j$ such
that this is true, the aggregation of $x_{mid(i)}+x_{mid(j)}$ adds
to flux across $x_{r}$. Therefore the discretized flux contributed
to $J_{r}^{N}(t_{k})$ by $x_{mid(i)}$ and a given $x_{mid(j)}$
is
\[
\Delta xg_{i}^{N}(t_{k})\int_{x_{j}}^{x_{j+1}}\frac{K_{A}(x_{mid(i)},y)}{y}g_{j}^{N}(t_{k})dy,
\]
and we can simply sum across all $j$ for which $x_{mid(j)}\ge\tilde{x}$.
The small exception to this rule occurs for the lowest $j$. In that
case, the lower limit of integration is $x_{mid(j)}$ instead of $x_{j}$.

Now we can use the FLFM to solve the Smoluchowski equation. For consistency,
$N$ represents the number of element boundaries, where each element
boundary is denoted $x_{i}$ for $i=1,\dots,N$ such that $x_{1}=x_{min}$
and $x_{N}=x_{max}$. Again we only consider uniform grids, so $\Delta x=\Delta x_{i}=x_{i+1}-x_{i\,}\forall i$.
Then at each time step and at each element midpoint, we approximate
(\ref{eq:smolu_flux_g}) with 
\begin{equation}
\frac{g_{i}^{N}(t_{k+1})-g_{i}^{N}(t_{k})}{\Delta t}=\frac{J_{i+1}^{N}(t_{k})-J_{i}^{N}(t_{k})}{\Delta x},\, i\in[1,N-1]\label{eq:flux deriv approx}
\end{equation}
and we approximate our initial conditions with
\begin{equation}
g_{i}^{N}(0)=\frac{1}{\Delta x}\int_{x_{i}}^{x_{i+1}}xf(0,x)dx.\label{eq:init cond flux}
\end{equation}
Now recall that $J[g]$ is defined as $\int_{0}^{x}\int_{x-u}^{x_{max}}\frac{K_{A}(u,v)}{v}g(t,u)g(t,v)dvdu$
with $x\in(0,x_{max})$ and $t\in\mathbb{R}_{+}$, so evaluation of
the right hand side of (\ref{eq:flux deriv approx}) is less obvious.
To illustrate how we calculate the right hand side of our second order
approximation in (\ref{eq:flux deriv approx}), we present it in Appendix
\ref{sec:flux calcs-1}. In \citep{Filbet2004}, Filbet and Laurençot
demonstrate second order convergence in $\Delta x$ of the FLFM in
the $L^{1}$ norm. Using similar arguments as those made in Section
\ref{sub:Num Meth FEM}, we note the convergence of the $i^{th}$
moments is second order (or faster) in $\Delta x$. 

\begin{rmk}
\it{The i$^{th}$ moment $M_i^N$ of the approximate solution converges to the i$^{th}$ moment $M_i$ of the analytical solution with order 2 (or faster) in $L^1[\bold{X}]$.}
\end{rmk}
\begin{proof}
Recall that for fixed $t_{k}$, there is second order convergence
in $\Delta x$ of $g^{N}$ to g.
\begin{eqnarray*}
\lim_{N\rightarrow\infty}\left|M_{i}^{N}(f(t_{k},\cdot);\mathbf{X})-M_{i}(f(t_{k},\cdot);\mathbf{X})\right| & = & \lim_{N\rightarrow\infty}\left|\int_{\mathbf{X}}x^{i}f^{N}(t_{k},x)dx-\int_{\mathbf{X}}x^{i}f(t_{k},x)dx\right|\\
 & = & \lim_{N\rightarrow\infty}\left|\int_{\mathbf{X}}x^{i-1}g^{N}(t_{k},x)dx-\int_{\mathbf{X}}x^{i-1}g(t_{k},x)dx\right|\\
 & = & \lim_{N\rightarrow\infty}\left|\int_{\mathbf{X}}x^{i-1}\left\{ g^{N}(t_{k},x)-g(t_{k},x)\right\} dx\right|\\
 & \le & \lim_{N\rightarrow\infty}\int_{\mathbf{X}}x^{i-1}\left|g^{N}(t_{k},x)-g(t_{k},x)\right|dx\\
 & \le & \int_{\mathbf{X}}x^{i-1}\lim_{N\rightarrow\infty}\left\Vert g^{N}(t_{k},x)-g(t_{k},x)\right\Vert _{L^{1}[x_{j},x_{j+1}]}
\end{eqnarray*}

\end{proof}

\section{Analytical Solutions\label{sec:Analytical-Solutions}}

In Figure 2 of \citep{Wattis2006a}, Wattis provides a diagram partitioning
regions of different generic behavior for varying aggregation kernels
with general form $K_{A}(x,y)=x^{\mu}y^{\nu}+x^{\nu}y^{\mu}$ where
$\mu,\nu\in\mathbb{\mathbb{R}}^{+}$, including the exactly solvable
cases. Lee generates similar conclusions with respect to the generic
behavior for varying aggregation kernels in \citep{Lee2001jopa}.
In order to make fair comparisons between the FEM and the FLFM, we
use known solutions to (\ref{eq:smolu}) for the two aggregation kernels,
$K_{A}(x,y)\equiv1$ and $K_{A}(x,y)=xy$. The aggregation kernel,
$K_{A}(x,y)\equiv1$, represents a system with slower aggregation
and no gelation, while $K_{A}(x,y)=xy$ represents a system with rapid
aggregation where gelation does occur. By including both kernels,
not only can we compare the FEM and FLFM to known solutions, but we
cover a breadth of possible systems. We also recognize that the information
lost in the interval $[0,x_{1})$ degrades the overall accuracy of
the respective methods, and we address this in more detail in Section
\ref{sub:Validation}. 

\begin{figure}[H]
\noindent \raggedright{}\subfloat[$K_{A}(x,y)\equiv1$ ]{\label{a is 1 analytic}

\begin{centering}
\includegraphics[scale=0.43]{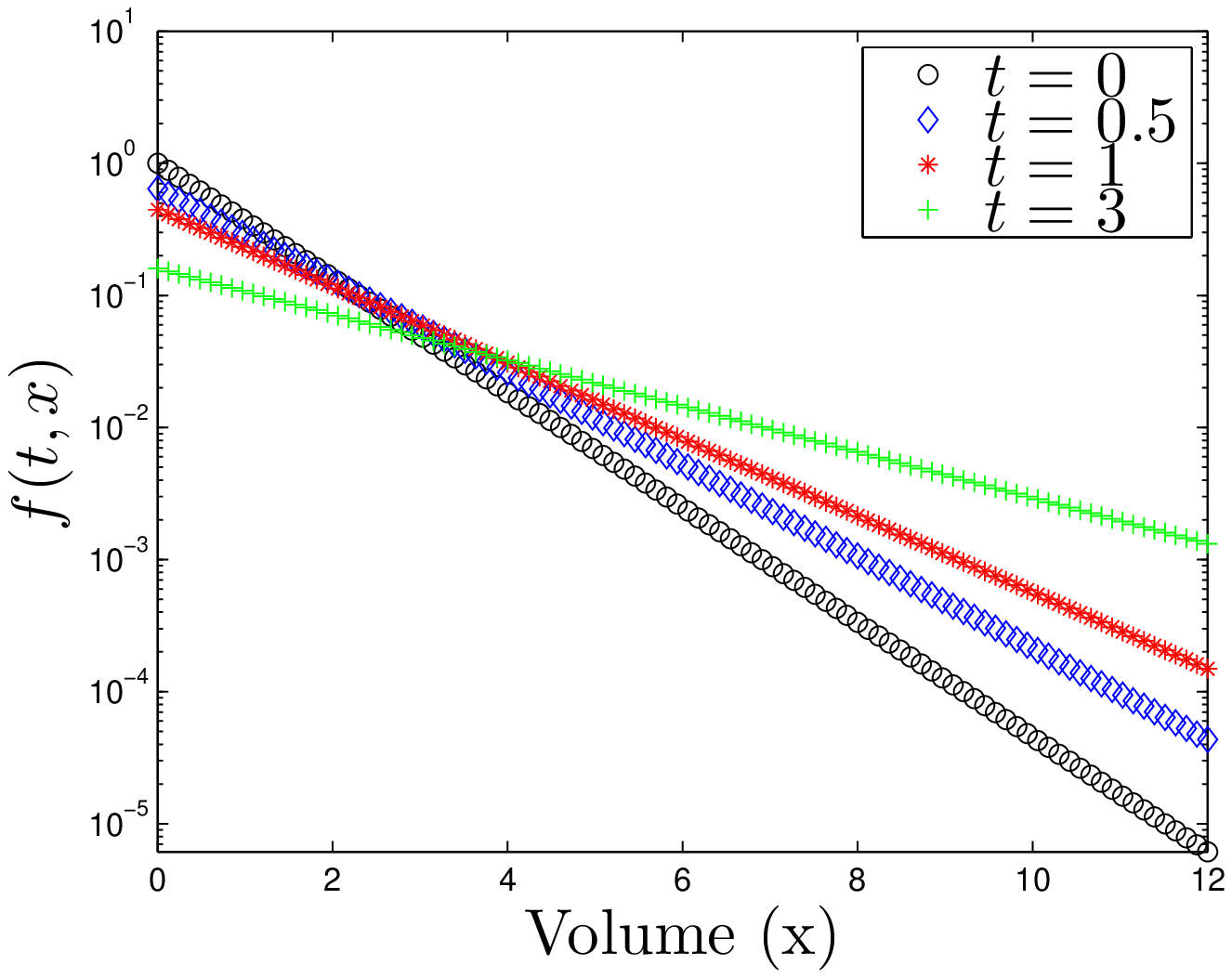}
\par\end{centering}

}\hfill{}\subfloat[$K_{A}(x,y)=xy$]{\label{a is xx' analytic}

\begin{centering}
\includegraphics[scale=0.43]{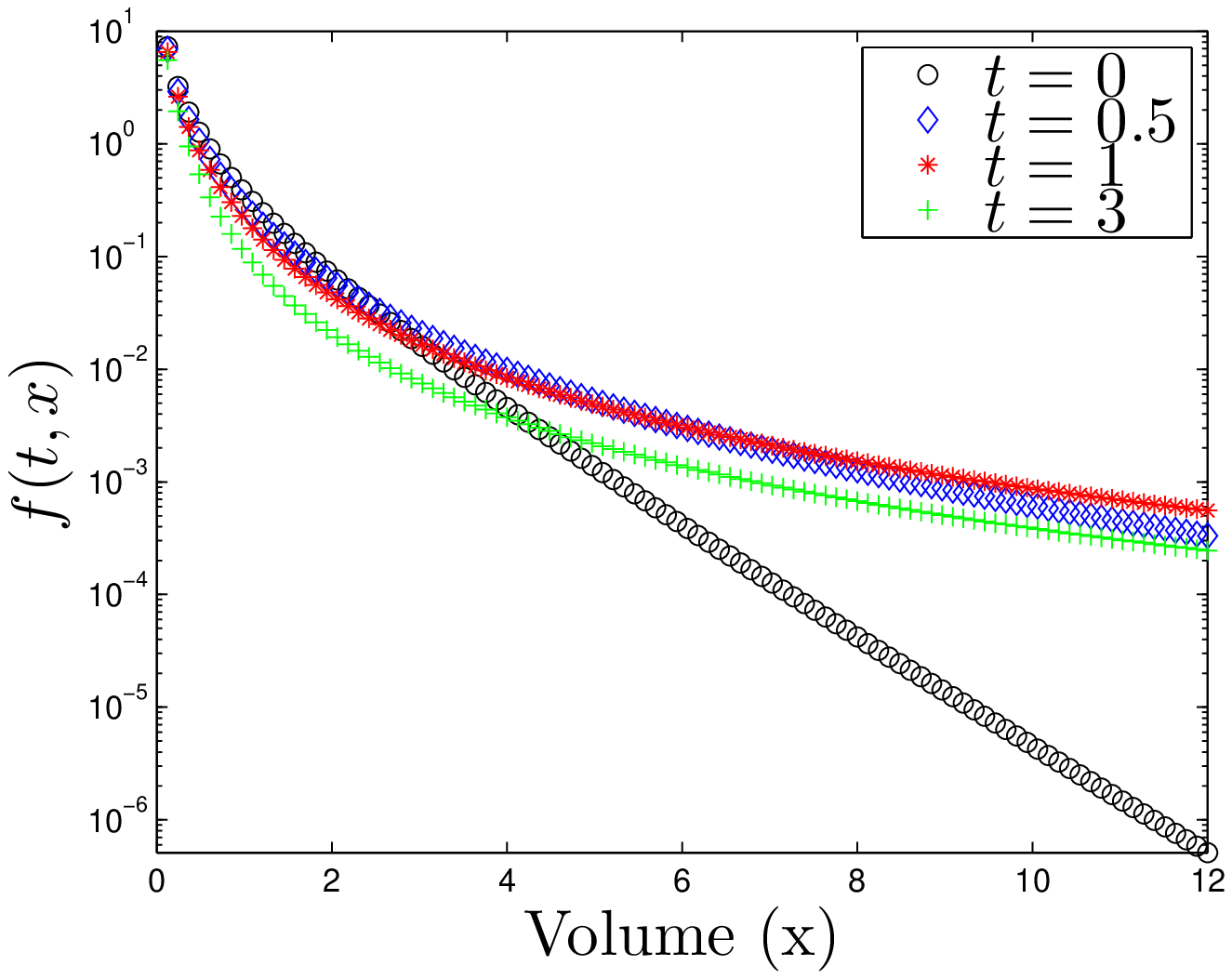}
\par\end{centering}

}\caption{Analytic solutions - size distribution, $f(t,x)$, vs volume, $x$,
for several snapshots in time: $t=0$, $t=0.5$, $t=1$, and $t=3$
in semilog scale. We use $K_{A}(x,y)\equiv1$ and $K_{A}(x,y)=xy$
as the two aggregation kernels representing a breadth of modeled systems
from which we can base comparisons of the two methods studied, the
FEM and the FLFM. }
\end{figure}

A known analytical solution for $K_{A}(x,y)\equiv1$ is \citep{Filbet2004}
\begin{equation}
f(t,x)=\left(\frac{2}{2+t}\right)^{2}e^{-\frac{2}{2+t}x}.\label{eq: truef a is 1}
\end{equation}
 We depict the solution (\ref{eq: truef a is 1}) in Figure \ref{a is 1 analytic}
for several different snapshots in time. Similarly, a known analytical
solution for $K_{A}(x,y)=xy$ is \citep{Ernst1984,Filbet2004} 
\begin{equation}
f(t,x)=e^{(-T(t)x)}\frac{I_{1}(2x\sqrt{t})}{x^{2}\sqrt{t}}.\label{eq:xx' solution}
\end{equation}
Here
\[
T(t)=\left\{ \begin{array}{ll}
1+t & \mbox{if}\,\, t\le1\\
2\sqrt{t} & \mbox{otherwise}
\end{array}\right.,
\]
and
\[
I_{1}(x)=\frac{1}{\pi}\int_{0}^{\pi}e^{x\,\cos\theta}\cos\theta d\theta
\]
is a modified Bessel function of the first kind. We depict in Figure
\ref{a is xx' analytic} the solution (\ref{eq:xx' solution}) for
several different snapshots in time. For this solution, note that
$f(0,x)=\frac{e^{-x}}{x}$, which is not necessarily obvious (see
Appendix \ref{sec:t=00003D0 solution notes} for the derivation). 

\section{Computational Results\label{sec:Computational-Results}}

Using the analytical solutions described in the previous section,
we can test both numeric schemes, the FEM and the FLFM, to compare
their results. In Section \ref{sub:Validation}, we define our measure
of error and compare the resulting convergence rates achieved by both
methods. In Section \ref{sub:Moment calcs}, we compare the FEM's
and the FLFM's accuracy in approximating the zeroth and first moments.
In Section \ref{sub:Compcost}, we discuss the computation cost required
for each of our simulations. Finally, in Section \ref{sub:Rel between xmax delx},
we examine the effects of grid spacing and our truncation parameter,
$x_{max}$.

\subsection{\label{sub:Validation}Validation}

\begin{figure}[H]
\subfloat[FEM]{\label{FEM error ais1}

\begin{centering}
\includegraphics[scale=0.44]{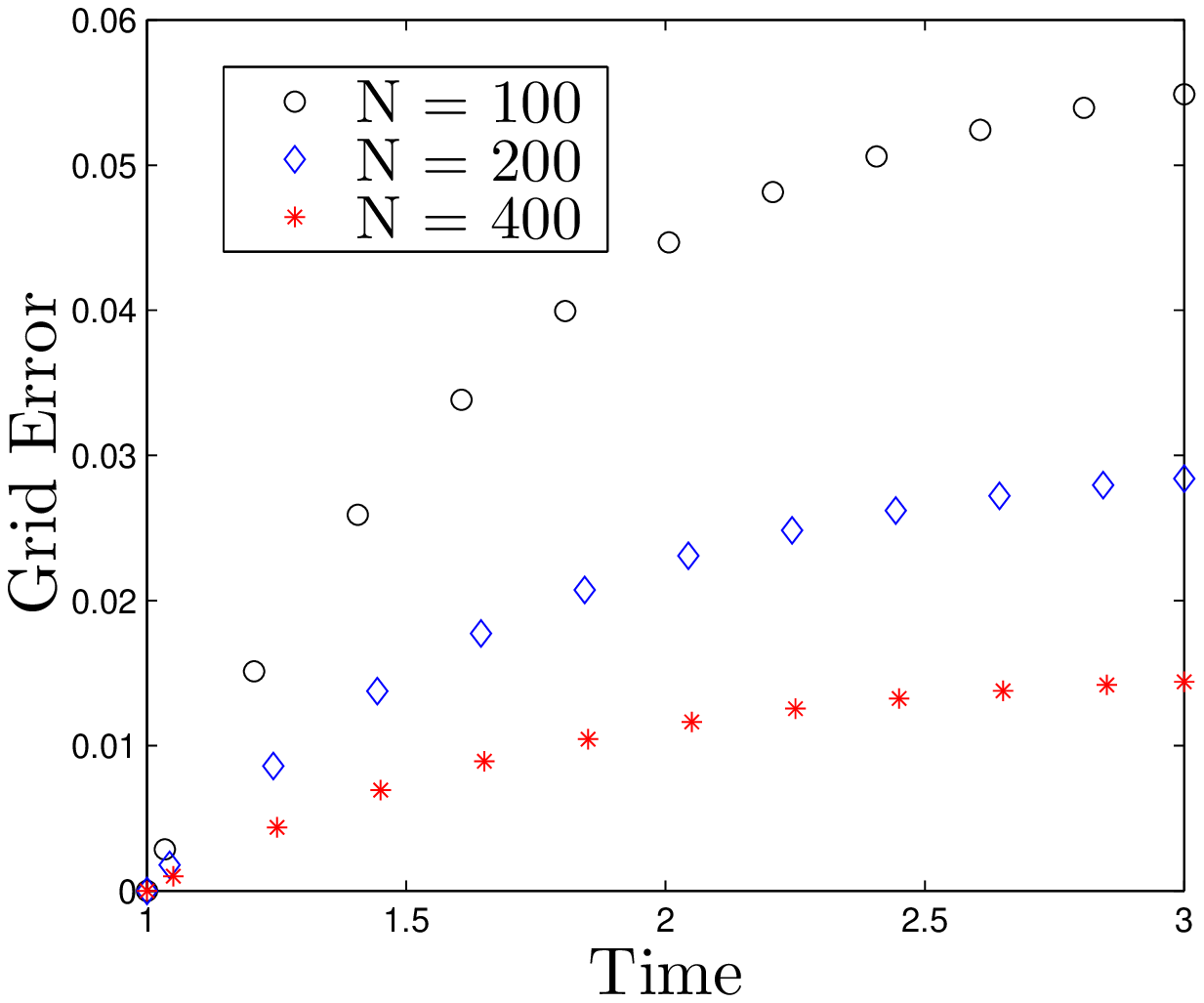}
\par\end{centering}

}\hfill{}\subfloat[FLFM]{\label{flux error ais1}

\begin{centering}
\includegraphics[scale=0.44]{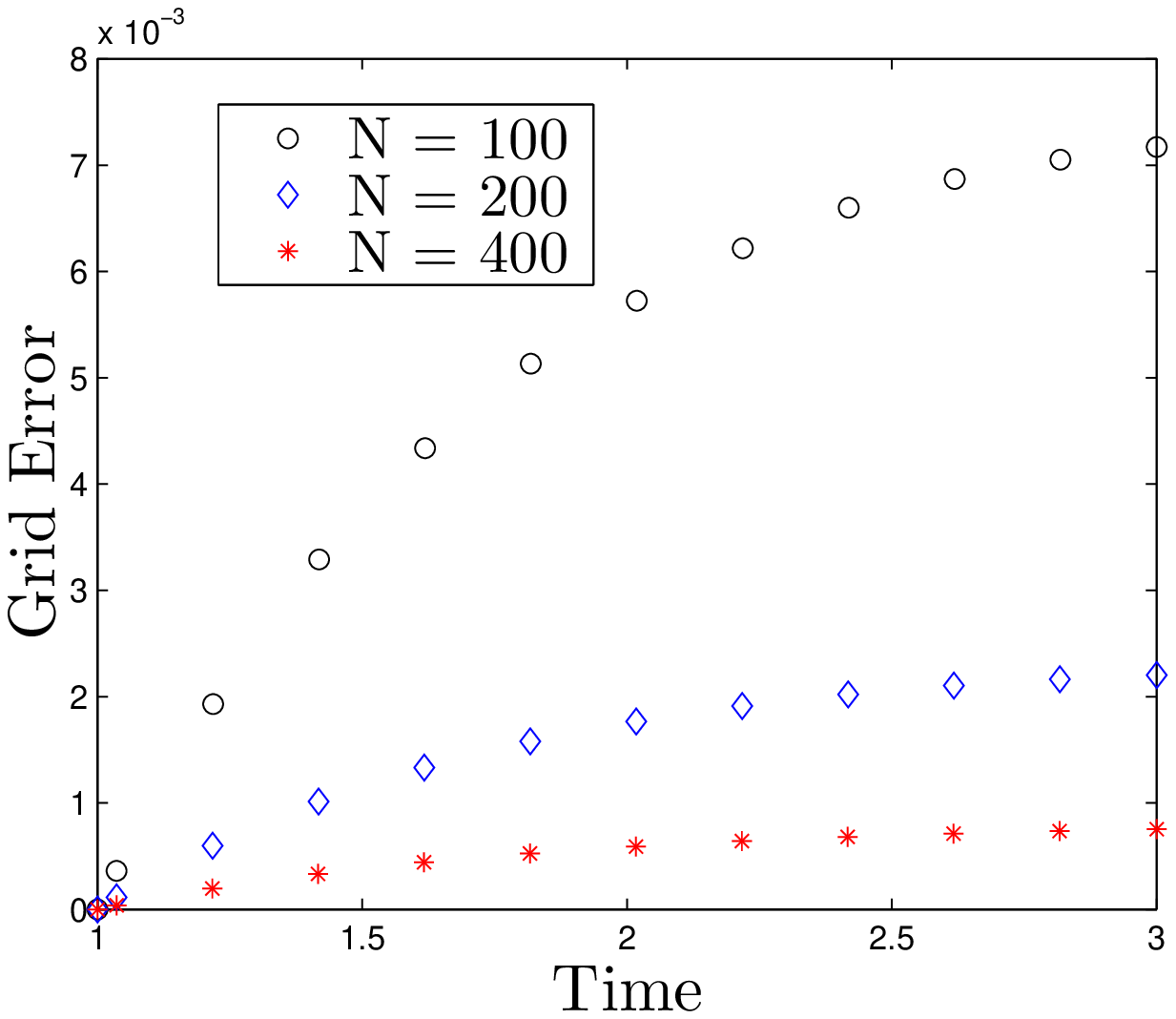}
\par\end{centering}

}\caption{Error vs.~Time - $L^{1}$ grid error norm vs.~time when $K_{A}(x,y)\equiv1$
for increasing grid densities: $N=100$, $N=200$, and $N=400$ in
linear scale. We define the error for the FEM based on the difference
between analytic and numeric size distribution, $f$, but we base
the error for the FLFM on the difference between analytic and numeric
volume distribution, $xf$. The FLFM provides a higher order of accuracy
with this aggregation kernel, which represents a system that aggregates
slowly and does not experience gelation.}
\label{error ais1}
\end{figure}

\begin{figure}[H]
\subfloat[FEM]{\label{FEM error aisxy}

\begin{centering}
\includegraphics[scale=0.44]{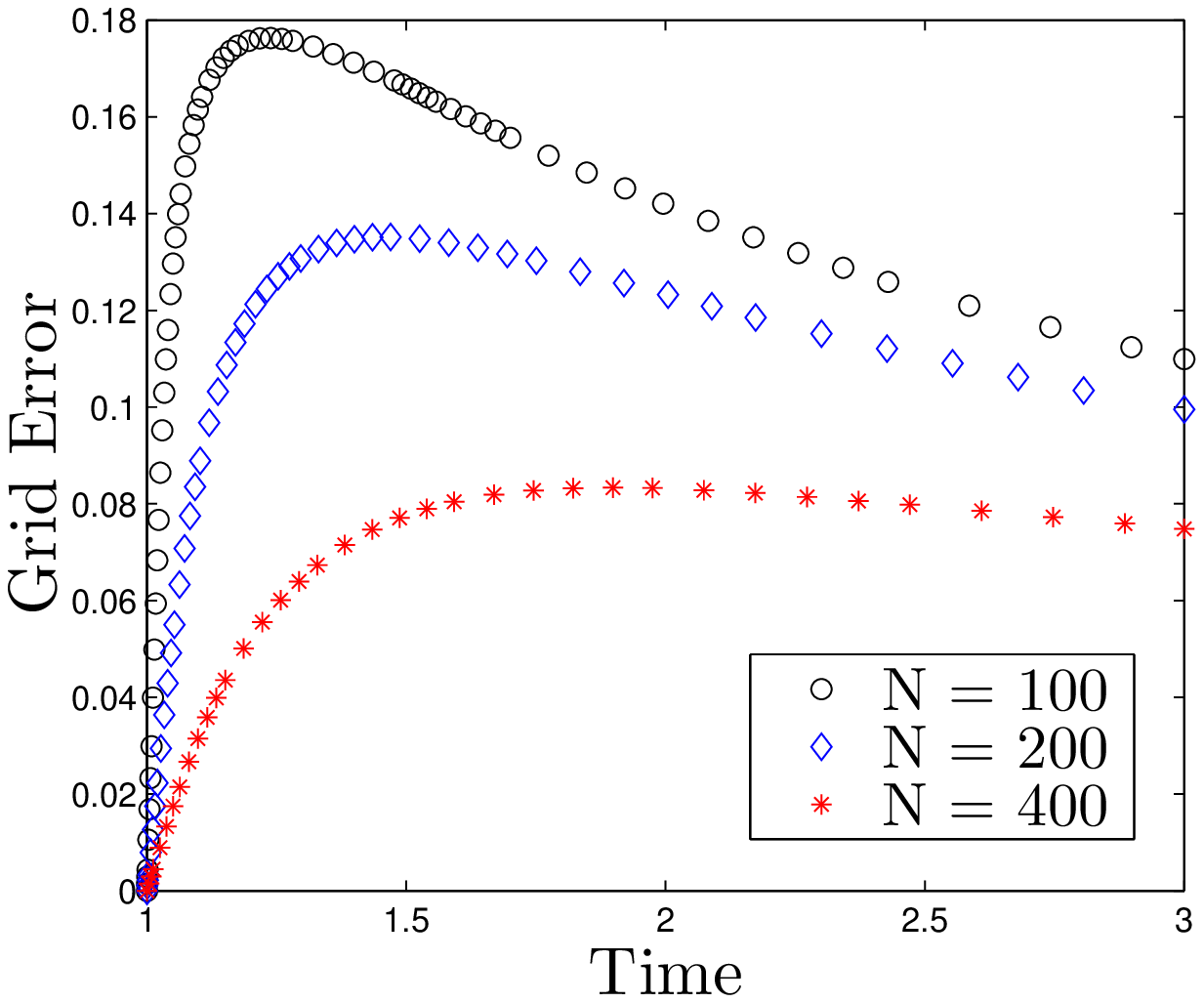}
\par\end{centering}

}\hfill{}\subfloat[FLFM]{\label{flux error aisxy}

\begin{centering}
\includegraphics[scale=0.44]{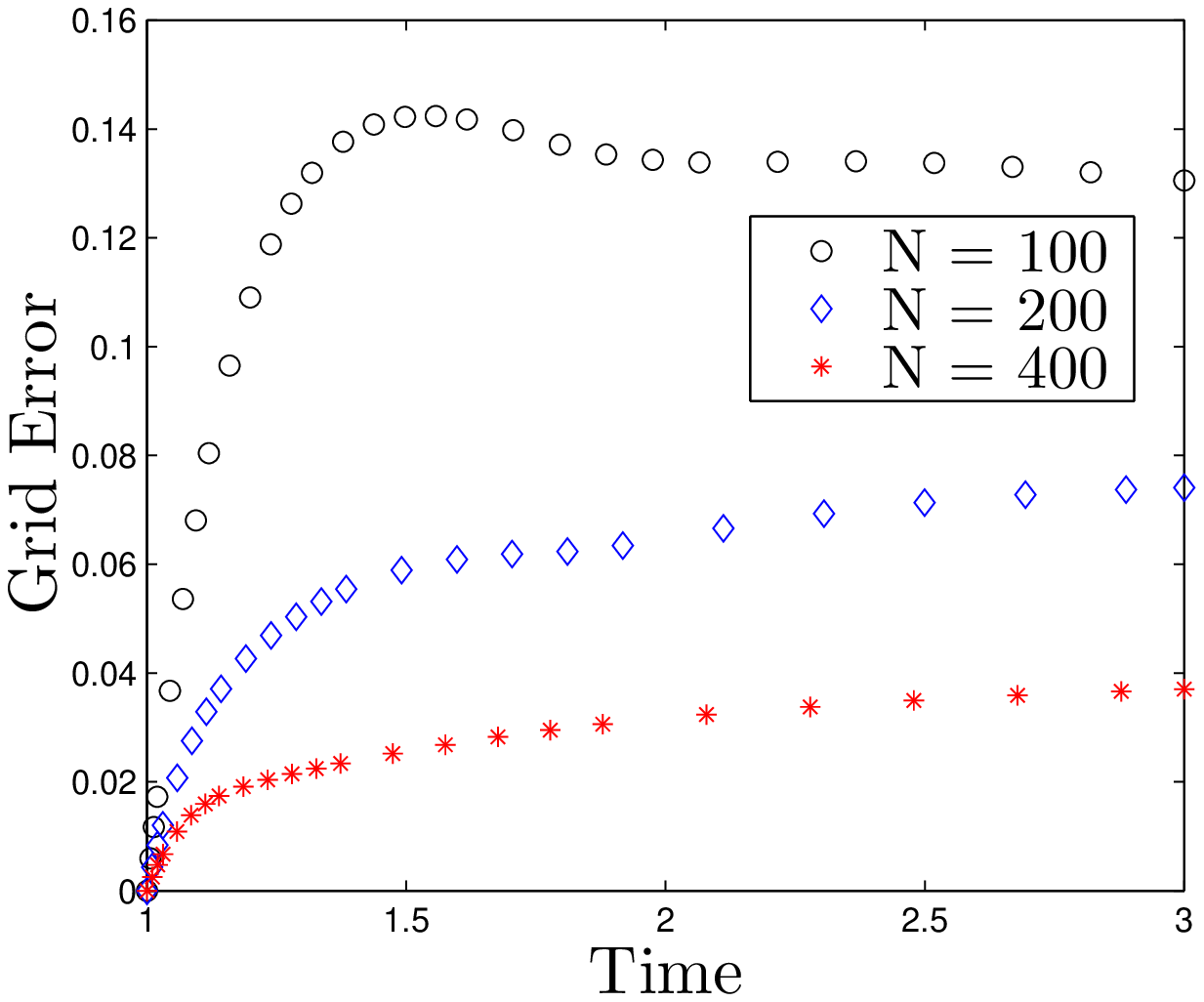}
\par\end{centering}

}\caption{Error vs.~Time - $L^{1}$ grid error norm vs.~time when $K_{A}(x,y)=xy$
for increasing grid densities: $N=100$, $N=200$, and $N=400$ in
linear scale. We define the error for the FEM based on the difference
between analytic and numeric size distribution, $f$, but we base
the error for the FLFM on the difference between analytic and numeric
volume distribution, $xf$. The FLFM provides a higher order of accuracy
than FEM as well for this aggregation kernel, which represents a rapidly
aggregating system where gelation occurs.}
\label{error aisxy}
\end{figure}

We cannot overstate the importance of computing the correct norm when
determining the error in a particular numerical scheme. In particular,
when using the FEM, we approximate $f$, but using the FLFM, we approximate
$g$. We then denote the mean analytical solution across the $i^{th}$
element at $t_{k}$ as $f_{i}(t_{k})$ and $g_{i}(t_{k})$ respectively.
To compute the error for both methods, we use the grid function norm%
\footnote{See Appendix A.5 of LeVeque \citep{LeVeque2007} for further details.%
}. For the FEM, we denote the error
\[
e_{i}(t_{k})=\left|f_{i}^{N}(t_{k})-\frac{M_{0}(f(t_{k},\cdot);x_{i},x_{i+1})}{\Delta x}\right|=\left|f_{i}^{N}(t_{k})-\pi^{N}(f(t_{k}))\right|
\]
and for the FLFM, we denote the error 
\[
e_{i}(t_{k})=\left|g_{i}^{N}(t_{k})-\frac{M_{1}(\frac{g(t_{k},\cdot)}{\cdot};x_{i},x_{i+1})}{\Delta x}\right|=\left|g_{i}^{N}(t_{k})-g_{i}(t_{k})\right|.
\]
In both cases, we generate vectors discretizing an error function,
$e(t_{k},x)$, so we calculate the grid function norm as
\[
\left\Vert e(t_{k})\right\Vert _{1}=\Delta x\sum_{i=1}^{N-1}e_{i}(t_{k}).
\]
With our error defined, we aim to reduce the impact of information
lost from excluding the interval $[0,x_{1})$ by using the domain
$t=[1,3]$. We now compare the accuracy of the two methods for each
aggregation kernel, and we plot the results in Figures \ref{error ais1}
and \ref{error aisxy}.

\begin{figure}[H]
\centering{}\includegraphics[scale=0.45]{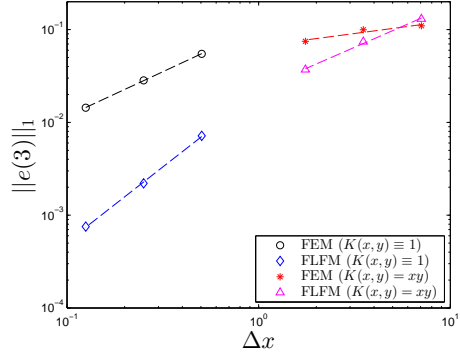}\caption{Error vs.~$\Delta x$ - $L^{1}$ grid error norm at $t=3$ vs.~$\Delta x$
in log scale for both the FEM and the FLFM using both aggregation
kernels, $K_{A}(x,y)\equiv1$ and $K_{A}(x,y)=xy$. The slope of the
dashed lines represent the approximate order of accuracy of each method
using each respective aggregation kernel. The FEM achieves approximately
first and 0.3 order accuracy for the kernels $K_{A}(x,y)\equiv1$
and $K_{A}(x,y)=xy$ respectively. The FLFM achieves approximately
1.5 and first order for the kernels $K_{A}(x,y)\equiv1$ and $K_{A}(x,y)=xy$
respectively.}
\label{ordersummary}
\end{figure}

When $K_{A}(x,y)\equiv1$, the average analytical solution across
the $i^{th}$ element is 
\begin{eqnarray*}
f_{i}(t_{k}) & = & \frac{1}{\Delta x}\int_{x_{i}}^{x_{i+1}}\left(\frac{2}{2+t}\right)^{2}e^{-\left(\frac{2}{2+t}\right)y}dy\\
 & = & \frac{1}{\Delta x}\left(-\frac{2}{2+t}\right)\left[e^{-\frac{2x_{i+1}}{2+t}}-e^{-\frac{2x_{i}}{2+t}}\right],
\end{eqnarray*}
from which we calculate the error in the FEM. Conversely, when we
calculate the error in the FLFM, the analytical solution across the
$i^{th}$ element is
\begin{eqnarray*}
g_{i}(t_{k}) & = & \frac{1}{\Delta x}\int_{x_{i}}^{x_{i+1}}\left(\frac{2}{2+t}\right)^{2}ye^{-\left(\frac{2}{2+t}\right)y}dy\\
 & = & \frac{1}{\Delta x}\left[\left(\frac{2x_{i}}{2+t}+1\right)\left(e^{-\frac{2x_{i}}{2+t}}\right)-\left(\frac{2x_{i+1}}{2+t}+1\right)\left(e^{-\frac{2x_{i+1}}{2+t}}\right)\right].
\end{eqnarray*}
With $K_{A}(x,y)\equiv1$, $x_{1}$ can be as small as we wish, but
physically, $x_{1}>0$, so we choose $x_{1}=10^{-3}$. We should also
note that for any choice of $x_{1}>0$, we introduce error by disregarding
information generated in the true solution by volumes in the range
$[0,x_{1}).$ Under these conditions, we achieve approximately first
order accuracy using the finite element approach as depicted in Figure
\ref{ordersummary}. We achieve approximately 1.5 order accuracy using
the flux approach, also depicted in Figure \ref{ordersummary}, and
the overall error is much smaller.

When $K_{A}(x,y)=xy$, to calculate error in the FEM approach, we
use the average analytical solution across the $i^{th}$ element
\begin{equation}
f_{i}(t_{k})=\frac{1}{\Delta x}\int_{x_{i}}^{x_{i+1}}e^{(-Ty)}\frac{I_{1}(2y\sqrt{t})}{y^{2}\sqrt{t}}dy.\label{eq:avg_f}
\end{equation}
To calculate the error in the FLFM approach, we use the average analytical
solution across the $i^{th}$ element 
\begin{equation}
g_{i}(t_{k})=\frac{1}{\Delta x}\int_{x_{i}}^{x_{i+1}}e^{(-Ty)}\frac{I_{1}(2y\sqrt{t})}{y\sqrt{t}}dy.\label{eq:avg_g}
\end{equation}
Neither (\ref{eq:avg_f}) or (\ref{eq:avg_g}) have analytic integrals,
so we apply global adaptive quadrature, as implemented in Matlab's
\textbf{integral2}, with the default relative tolerance of $10^{-6}$
and absolute tolerance of $10^{-10}$, to approximate $f_{i}$ and
$g_{i}$ respectively at each time step. Additionally, extremely small
volume sizes create large numerical inaccuracies in these integrals,
so to achieve the best results, we choose $x_{1}=0.75$. Under these
assumptions, we achieve approximately 0.3 order accuracy using the
FEM approach, whereas we achieve approximately first order accuracy
using FLFM as depicted in Figure \ref{ordersummary}. For both methods,
and with both aggregation kernels, we achieve less accuracy than the
maximum achievable by the respective methods, which we address below.
 
\begin{figure}[H]
\centering{}\includegraphics[scale=0.45]{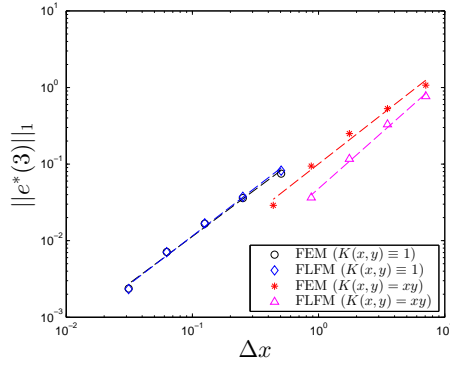}\caption{Error vs.~$\Delta x$ - $L^{1}$ grid error norm at $t=3$ vs.~$\Delta x$
in log scale for both the FEM and the FLFM using the aggregation kernel,
$K_{A}(x,y)=xy$. The slope of the dashed lines represents the approximate
order of accuracy of each method. Both methods trend towards approximately
second order accuracy.}
\label{800ordersummary}
\end{figure}

Undoubtedly, the information lost in the interval $[0,x_{1})$ degrades
the overall accuracy of the respective methods. In light of these
inaccuracies, we offer another test of the methods' convergence rates
for both aggregation kernels. In this second test, we compare solutions
for both kernels using 100, 200, 400, 800, and 1600 points to a fine
grid solution of 3200 points.%
\footnote{Due to the large computation time required by the FLFM when $K(x,y)=xy,$
we currently present that case's results for 100, 200, 400, and 800
grid points compared to a fine grid solution of 1600 points in Figure
\ref{800ordersummary}.%
} The error here is analogous to the error measured when we compare
our approximated solutions for varying $N$ to the analytic solution.
In the interest of clarity, we use a superscript, {*}, to denote the
error measured when we compare our approximated solutions for varying
$N$ to the fine grid solution where $N=3200.$ For the FEM, we then
denote the error
\[
e_{i}^{*}(t_{k})=\left|f_{i}^{N}(t_{k})-f_{i}^{3200}(t_{k})\right|,
\]
and for the FLFM, we denote the error 
\[
e_{i}^{*}(t_{k})=\left|g_{i}^{N}(t_{k})-g_{i}^{3200}(t_{k})\right|.
\]
In both cases, we generate vectors discretizing an error function,
$e^{*}(t_{k},x)$, so we calculate the grid function norm as
\[
\left\Vert e^{*}(t_{k})\right\Vert _{1}=\Delta x\sum_{i=1}^{N-1}e_{i}^{*}(t_{k}).
\]
With this test, we achieve an order of accuracy that trends towards
second order with each doubling of the number of grid points. These
results support our conjecture of second order accuracy for the FEM
and support the work in \citep{Filbet2004} that demonstrates second
order accuracy for the FLFM. We depict these results in Figure \ref{800ordersummary}.

\subsection{\label{sub:Moment calcs}Moment calculations}

In this section, we investigate both methods' accuracy in approximating
the moments. Our study reveals that an inherent advantage of the FEM
in its use of a size distribution results in a much more accurate
approximation of the zeroth moment when $K_{A}(x,y)\equiv1$. Conversely,
the FLFM's use of a volume distribution results in a more accurate
approximation of the first moment for both aggregation kernels.

Mathematically, we represent total aggregates as the zeroth moment,
$M_{0}$, and total volume as the first moment, $M_{1}$, as defined
in (\ref{eq:ith moment}). Different formulations of the governing
Smoluchowski equation induce different approaches to calculating the
moments numerically. For example, Guy, Fogelson, and Keener \citep{Guy2007}
extend the generating function approach described in \citep{Ziff1980}
to study blood clots. In the FEM approach, we approximate the size
distribution across the $i^{th}$ element as $f_{i}^{N}(t_{k})$,
whereas with the FLFM, we approximate the volume distribution across
the $i^{th}$ element as $g_{i}^{N}(t_{k})$. Using (\ref{eq:ith moment})
with our known analytical solutions for $f(t,x)$, the analytical
total number of particles in our truncated system is
\begin{equation}
M_{0}(t_{k})=\int_{x_{1}}^{x_{N}}f(t_{k},y)dy=\int_{x_{1}}^{x_{N}}\frac{g(t_{k},y)}{y}dy,\label{eq:M0redef}
\end{equation}
and the analytical total volume in our truncated system is
\begin{equation}
M_{1}(t_{k})=\int_{x_{1}}^{x_{N}}yf(t_{k},y)dy=\int_{x_{1}}^{x_{N}}g(t_{k},y)dy.\label{eq:M1redef}
\end{equation}
We acknowledge the mild abuse of the notation in redefining $M_{i}$
in (\ref{eq:M0redef}) and (\ref{eq:M1redef}) in the name of focusing
attention on how the moment changes as a function of time. In what
follows we will clearly note when referring to moments of the FEM
or the FLFM solutions. We demonstrate below that the FEM more accurately
approximates the zeroth moment when $K_{A}(x,y)\equiv1$. 

\begin{figure}[H]
\subfloat[$K_{A}(x,y)\equiv1$]{\label{M0ais1-1}

\centering{}\includegraphics[scale=0.45]{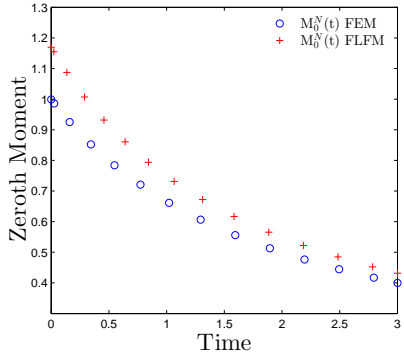}}\hfill{}\subfloat[$K_{A}(x,y)=xy$]{\label{M0aisxy-1}

\centering{}\includegraphics[scale=0.45]{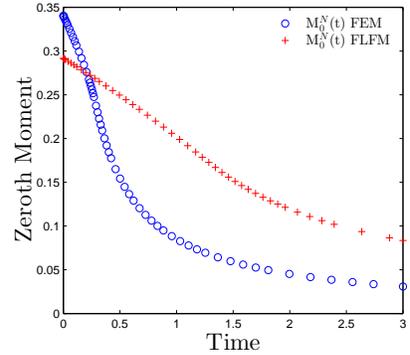}}\caption{$M_{0}^{400}(t)$ vs. Time - Numerical approximation of the zeroth
moment for a truncated volume domain, $x_{1}$ to $x_{400}$, for
both aggregation kernels. Note, with the FEM, $M_{0}^{400}(0)=M_{0}(0),$
but with the FLFM, $M_{0}^{400}(0)\approx M_{0}(0)$. When given data
in the form of a size distribution, the FEM's approximation of the
total particles in a system starts with an advantage over the FLFM's
approximation, which has error at $t=0$. We use a linear scale for
both plots.}
\label{M0}
\end{figure}

\begin{figure}[H]
\subfloat[Zeroth Moment]{\centering{}\includegraphics[scale=0.43]{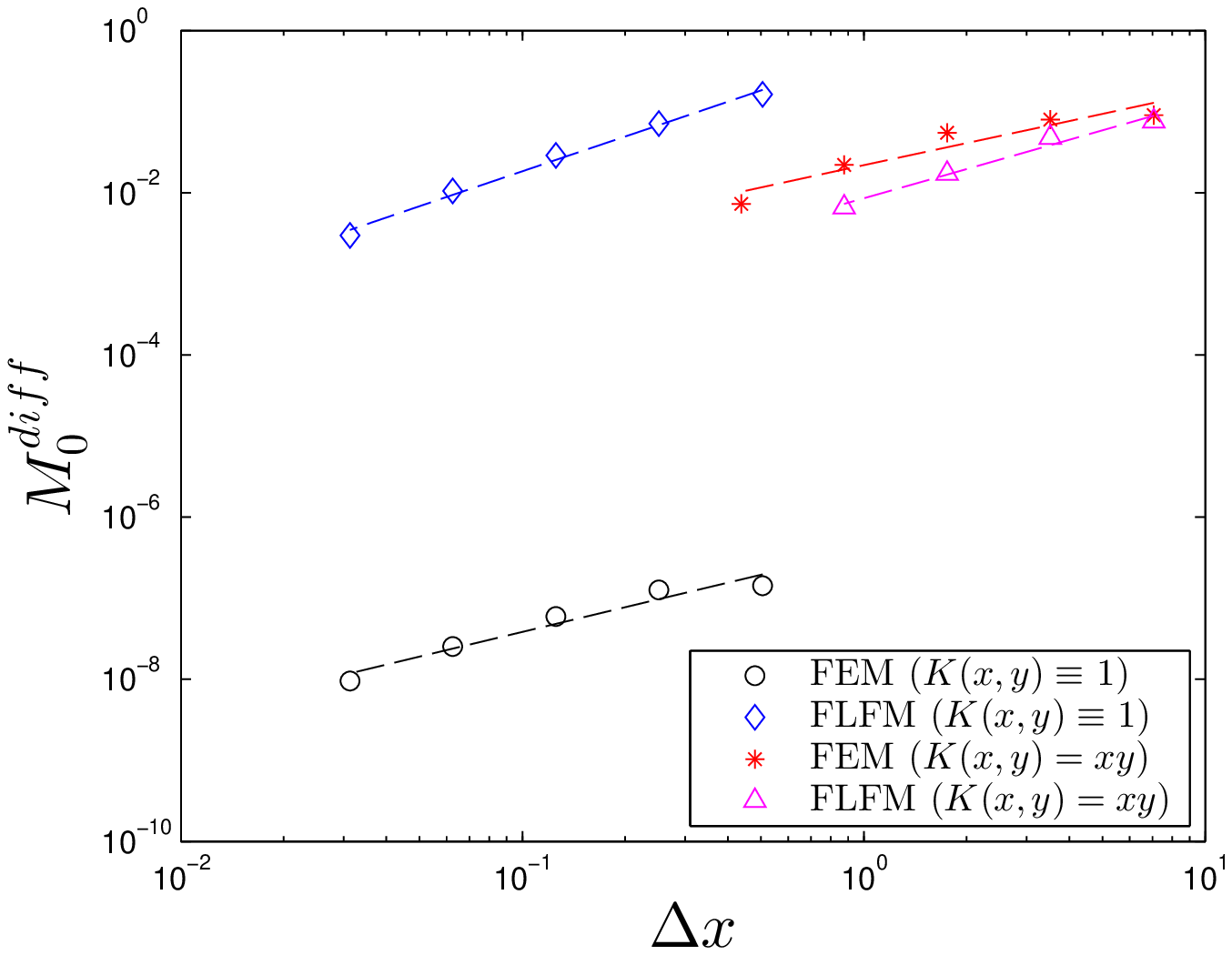}\label{M0ordsumm}}\hfill{}\subfloat[First Moment]{\centering{}\includegraphics[scale=0.43]{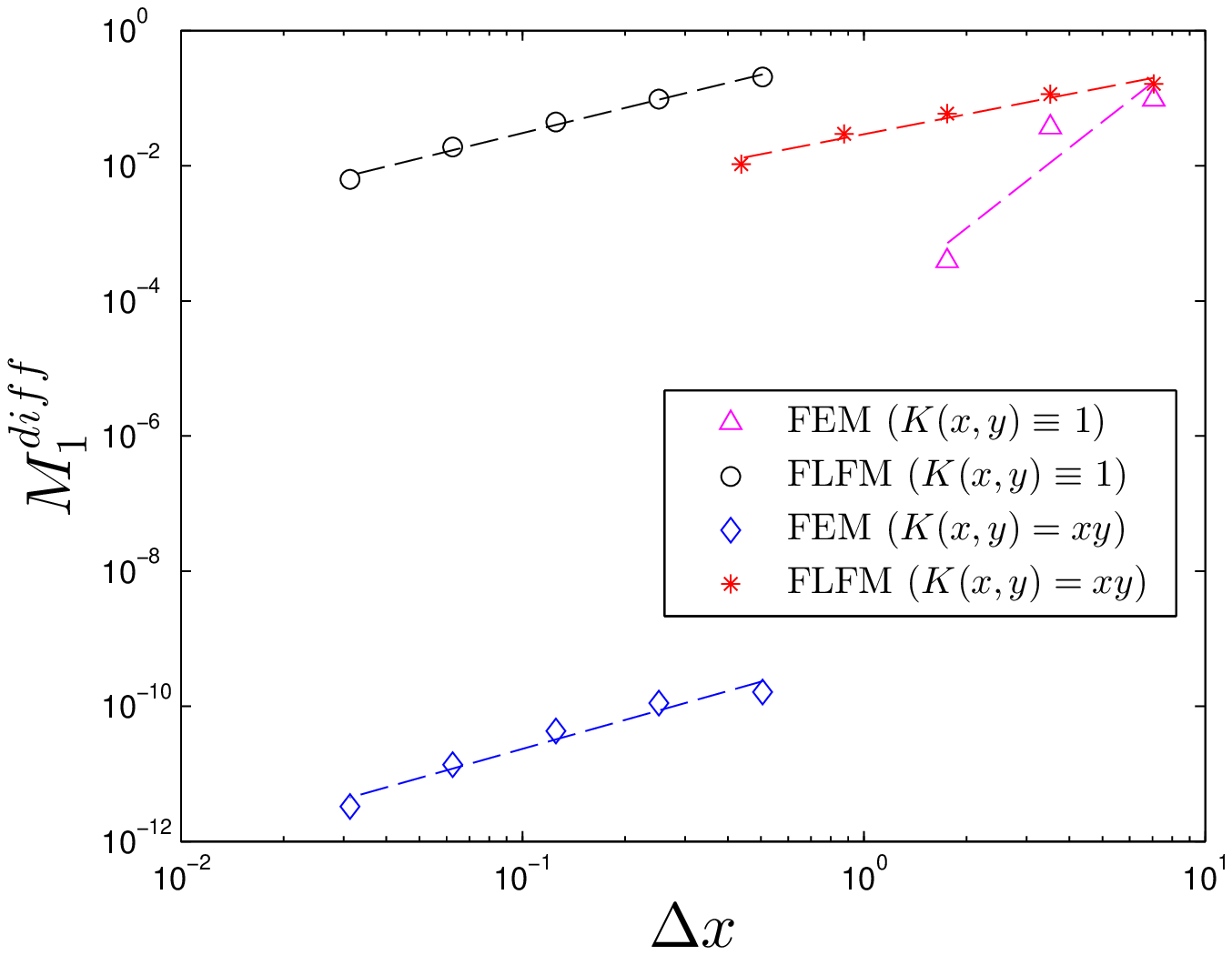}\label{M1ordsumm}}\caption{Moment differences at $t=3$ vs.~$\Delta x$ - Difference between
approximations of the zeroth and first moments by coarse grid ($N=100,\,200,\,400,\,800,\,1600$)
simulations and the approximation by the fine grid ($N=3200$) simulation.
In all cases, we observe a trend towards second order convergence,
with one case (slow aggregation, e.g., $K_{A}(x,y)\equiv1$ and zeroth
moment approximation) where the FEM is more accurate.}
\label{momsummpic}
\end{figure}

Through contrasting our approximations of the zeroth moment by the
two methods studied in this paper, we derive an important advantage
that the FEM possesses at $t=0$ for a slowly aggregating system (e.g.,
$K_{A}(x,y)\equiv1$). In this case, for the FEM,
\[
M_{0}^{N}(t_{k})=\sum_{i=1}^{N-1}\int_{x_{i}}^{x_{i+1}}e^{-y}dy=M_{0}(0),
\]
but for the FLFM
\[
M_{0}^{N}(0)=\sum_{i=1}^{N-1}\left[\left(\frac{1}{\Delta x}\left\{ \left(x_{i}+1\right)e^{-x_{i}}-\left(x_{i+1}+1\right)e^{-x_{i+1}}\right\} \right)\ln\left(\frac{x_{i+1}}{x_{i}}\right)\right]\approx M_{0}(0).
\]
To derive this result, we start with the general formulation of the
zeroth moment at $t_{k}$. For the FEM,
\[
M_{0}^{N}(t_{k})=\int_{x_{1}}^{x_{N}}f_{i}^{N}(t_{k})dy=\Delta x\sum_{i=1}^{N-1}f_{i}^{N}(t_{k}),
\]
whereas with the FLFM
\[
M_{0}^{N}(t_{k})=\int_{x_{1}}^{x_{N}}\frac{g_{i}^{N}(t_{k})}{y}dy=\sum_{i=1}^{N-1}\left[g_{i}^{N}(t_{k})\ln\left(\frac{x_{i+1}}{x_{i}}\right)\right].
\]
Then at $t=0$ and with $K_{A}(x,y)\equiv1$, we have $f(0,x)=e^{-x}$
and $g(0,x)=xe^{-x}$, which we initialize numerically via (\ref{eq:init cond FEM})
and (\ref{eq:init cond flux}) respectively. The analytical zeroth
moment,
\[
M_{0}(0)=\int_{x_{1}}^{x_{N}}f(0,y)dy=\sum_{i=1}^{N-1}\int_{x_{i}}^{x_{i+1}}e^{-y}dy,
\]
matches the approximation by the FEM,
\[
M_{0}^{N}(0)=\int_{x_{1}}^{x_{N}}f_{i}^{N}(0)dy=\Delta x\sum_{i=1}^{N-1}f_{i}^{N}(0)=\sum_{i=1}^{N-1}\int_{x_{i}}^{x_{i+1}}e^{-y}dy.
\]
Conversely, the approximation by the FLFM,
\begin{eqnarray*}
M_{0}^{N}(0) & = & \int_{x_{1}}^{x_{N}}\frac{g_{i}^{N}(0)}{y}dy=\sum_{i=1}^{N-1}\left[g_{i}^{N}(0)\ln\left(\frac{x_{i+1}}{x_{i}}\right)\right]\\
 & = & \sum_{i=1}^{N-1}\left[\left(\frac{1}{\Delta x}\left\{ \left(x_{i}+1\right)e^{-x_{i}}-\left(x_{i+1}+1\right)e^{-x_{i+1}}\right\} \right)\ln\left(\frac{x_{i+1}}{x_{i}}\right)\right],
\end{eqnarray*}
is clearly not the same as $M_{0}(0)$. Not surprisingly, our outputs
for $M_{0}^{N}(0)$ vary slightly in Figure \ref{M0ais1-1} such that
$M_{0}^{N}(0)=M_{0}(0)\approx0.999$ using the FEM, but $M_{0}^{N}(0)\approx1.169$
using the FLFM. Along similar lines, when $K_{A}(x,y)=xy$, we have
$M_{0}(0)$ matching the approximation by the FEM where $M_{0}^{N}(0)\approx0.34$,
but the FLFM produces $M_{0}^{N}(0)\approx0.2916$ as depicted in
Figure \ref{M0aisxy-1}. We summarize these results in Table \ref{ith moments time0}.
Clearly, with experimental data in the form of a size distribution,
the FLFM starts at a disadvantage approximating the zeroth moment
since its approximation at $t=0$ contains error. Furthermore, the
FEM maintains a better approximation of the zeroth moment as depicted
in Figure \ref{M0ordsumm}.

Now contrasting our approximations of the first moment, we demonstrate
more accuracy by the FLFM than by the FEM. Using the FEM, 
\[
M_{1}^{N}(t_{k})=\int_{x_{1}}^{x_{N}}yf(t_{k},y)dy\approx\frac{1}{2}\sum_{i=1}^{N-1}f_{i}^{N}(t_{k})(x_{i+1}^{2}-x_{i}^{2}),
\]
whereas using the FLFM,
\[
M_{1}^{N}(t_{k})\approx\int_{x_{1}}^{x_{N}}g_{i}^{N}(t_{k})dy=\Delta x\sum_{i=1}^{N-1}g_{i}^{N}(t_{k}).
\]
In this case, with $t=0$ and $K_{A}(x,y)\equiv1$, $M_{1}^{N}(0)=M_{1}(0)\approx1.0000$
using the FLFM, but $M_{1}^{N}(0)\approx1.0013$ using the FEM. Similarly,
when $t=0$ and $K_{A}(x,y)=xy$, $M_{1}(0)$ matches the approximation
using the FLFM where $M_{1}^{N}(0)\approx0.472$, but the FEM produces
$M_{1}^{N}(0)\approx0.602$. We summarize these results in Table \ref{ith moments time0}.
\begin{table}[H]
\subfloat[$K_{A}(x,y)\equiv1$]{\centering{}%
\begin{tabular}{cccc}
\hline 
 & Analytical & FEM & FLFM\tabularnewline
\hline 
$M_{0}(0)$ & 0.999 & 0.999 & 1.169\tabularnewline
$M_{1}(0)$ & 1.0000 & 1.0013 & 1.0000\tabularnewline
\end{tabular}\label{ith moment comp ais1}}\hfill{}\subfloat[$K_{A}(x,y)=xy$\label{ith moment comp aisxy-2}]{\centering{}%
\begin{tabular}{cccc}
\hline 
 & Analytical & FEM & FLFM\tabularnewline
\hline 
$M_{0}(0)$ & 0.34 & 0.34 & 0.2916\tabularnewline
$M_{1}(0)$ & 0.472 & 0.602 & 0.472\tabularnewline
\end{tabular}}\caption{Comparison of analytical partial $i^{th}$ moments ($i=0,1$) at $t=0$
with the approximations by both the FEM and the FLFM. The analytical
$M_{0}(0)$ exactly matches that of the FEM, but the FLFM's approximation
of $M_{0}(0)$ contains error. Conversely, the analytical $M_{1}(0)$
exactly matches that of the FLFM, but the FEM's approximation of $M_{1}(0)$
contains error.}
\label{ith moments time0}
\end{table}

Having illuminated the respective advantages the FEM has approximating
the zeroth moment, and the FLFM has approximating the first moment,
we now examine the convergence rates of the two methods to a fine
grid ($N=3200$) approximation of the moments. In this case, we compute
the difference between the lower resolution grid ($N=100,\,200,\,400,\,800,\,1600)$
approximations of the moments at $t=3$, and the fine grid ($N=3200$)
approximation.%
\footnote{Due to the large computation time required by the FLFM when $K(x,y)=xy,$
we currently present that case's results for 100, 200, 400, and 800
grid points compared to a fine grid solution of 1600 points in Table
\ref{Momconvrates} and Figure \ref{momsummpic}.%
} We denote the difference in the approximations of the zeroth moment,
$M_{0}^{diff}$, and the difference in the approximations of the first
moment, $M_{1}^{diff}$ where
\begin{eqnarray*}
M_{0}^{diff} & = & \left|M_{0}^{N}(3)-M_{0}^{3200}(3)\right|\\
M_{1}^{diff} & = & \left|M_{1}^{N}(3)-M_{1}^{3200}(3)\right|.
\end{eqnarray*}
Our simulations include both aggregation kernels, and we depict the
convergence rates in Figures \ref{M0ordsumm} and \ref{M1ordsumm}.
As expected, for both methods, we observe a trend towards second order
accuracy (more evidence supporting our claim of second order convergence
for the FEM). Intriguingly, when $K_{A}(x,y)\equiv1$, the FEM more
accurately predicts the zeroth moment. In this case, given experimental
data in the form of a size distribution and a system experiencing
slow aggregation, the FEM is a better choice of methods. We summarize
convergence rates of these simulations in Table \ref{Momconvrates},
and note that as the number of grid points double the convergence
rates tend toward the expected convergence rates described in Sections
\ref{sub:Num Meth FEM} and \ref{sub:Num Meth FLFM}.
\begin{table}[H]
\centering{}%
\begin{tabular}{ccccccc}
\hline 
 & Method & Moment & 100$\rightarrow$200 & 200$\rightarrow$400 & 400$\rightarrow$800 & 800$\rightarrow$1600\tabularnewline
\hline 
\multirow{4}{*}{$K\equiv1$} & \multirow{2}{*}{FEM} & Zeroth & 0.2 & 1.1 & 1.2 & 1.4\tabularnewline
 &  & First & 1.1 & 1.1 & 1.2 & 1.6\tabularnewline
 & \multirow{2}{*}{FLFM} & Zeroth & 1.2 & 1.3 & 1.5 & 1.8\tabularnewline
 &  & First & 0.6 & 1.4 & 1.7 & 2\tabularnewline
\hline 
\multirow{4}{*}{$K=xy$} & \multirow{2}{*}{FEM} & Zeroth & 0.2 & 0.5 & 1.3 & 1.6\tabularnewline
 &  & First & 0.5 & 1.0 & 1.0 & 1.5\tabularnewline
 & \multirow{2}{*}{FLFM} & Zeroth & 0.7 & 1.5 & 1.4 & {*}$^{4}$\tabularnewline
 &  & First & 1.4 & 6.6 & {*}\footnotemark{} & {*}$^{4}$\tabularnewline
\end{tabular}\caption{Convergence rates of coarse grid ($N=100,\,200,\,400,\,800,\,1600)$
approximations of the zeroth and first moments to a fine grid ($N=3200$)
approximation of the moments. In all cases, we observe a trend towards
second order convergence.}
\label{Momconvrates}
\end{table}
\footnotetext{Results not reported due to limitations of machine precision and available computational resources\label{fnlabel}}

\subsection{Computation cost\label{sub:Compcost}}

Clearly, order of accuracy is important, but we also want to know
which method requires more computation in terms of floating point
operations. To make the comparisons, we compute the floating point
operations for the simulations discussed in Section \ref{sub:Validation}
for both the FEM and the FLFM. The number of operations counted includes
only the computations required for each algorithm to solve the system
of ODEs at a given point in time, $t_{k}$. We implement the simulations
using Matlab version 2013a on an Intel(R) Core(TM) i-5 2410M CPU @
2.3 GHz. For each simulation, the FEM requires significantly less
floating point operations than the FLFM. The results are summarized
in Table \ref{computation cost} and give us important insights. When
$K_{A}(x,y)\equiv1$, the required number of floating point operations
and an extrapolation of the error data in Figure \ref{ordersummary}
imply that the FEM with 800 grid points can achieve nearly equal accuracy
as the FLFM with 100 grid points. The FEM achieves this accuracy with
only 38\% of the floating point operations that it takes the FLFM,
so when computation cost is more valuable to the user, the FEM provides
a better choice.
\begin{table}[H]
\centering{}%
\begin{tabular}{cccccc}
\hline 
 &  & $N=100$ & $N=200$ & $N=400$ & $N=800$\tabularnewline
\hline 
\multirow{2}{*}{$K(x,y)\equiv1$} & FEM & $1.46\times10^{4}$ & $5.93\times10^{4}$ & $2.39\times10^{5}$ & $9.57\times10^{5}$\tabularnewline
 & FLFM & $2.49\times10^{6}$ & $1.99\times10^{7}$ & $1.60\times10^{8}$ & $1.28\times10^{9}$\tabularnewline
\hline 
\multirow{2}{*}{$K(x,y)=xy$} & FEM & $6.71\times10^{4}$ & $2.74\times10^{5}$ & $1.11\times10^{6}$ & $4.46\times10^{6}$\tabularnewline
 & FLFM & $2.48\times10^{6}$ & $1.99\times10^{7}$ & $1.60\times10^{8}$ & $1.28\times10^{9}$\tabularnewline
\end{tabular}\caption{Comparison of floating point operations required to solve the right
hand side system of ODEs for any $t_{k}$ for both the FEM and the
FLFM. We perform the simulations on a 2.3 GHz processor. When $K(x,y)\equiv1$,
the number of floating point operations performed for the FEM simulation
with 800 grid points requires only 38\% of the number of floating
point operations that the FLFM simulation with 100 grid points requires.
The error is nearly the same for both simulations implying that if
computation cost is more valuable to the user, the FEM should be the
choice of methods.}
\label{computation cost}
\end{table}

\subsection{\label{sub:Rel between xmax delx}Relationship between $x_{max}$
and $\Delta x$}

In \citep{Filbet2004}, Filbet and Laurençot report sensitivity of
the FLFM to the truncation parameter, $x_{max}.$ We study the error
as a function of both $x_{max}$ and $\Delta x$, and in particular,
we would like determine to what extent the truncation parameter, $x_{max}$,
impacts the overall accuracy of the two schemes. We do not find the
following results conclusive, but they do reveal intriguing patterns
that we plan to study more extensively in future work.

\begin{figure}[H]
\subfloat[FEM]{\centering{}\includegraphics[scale=0.45]{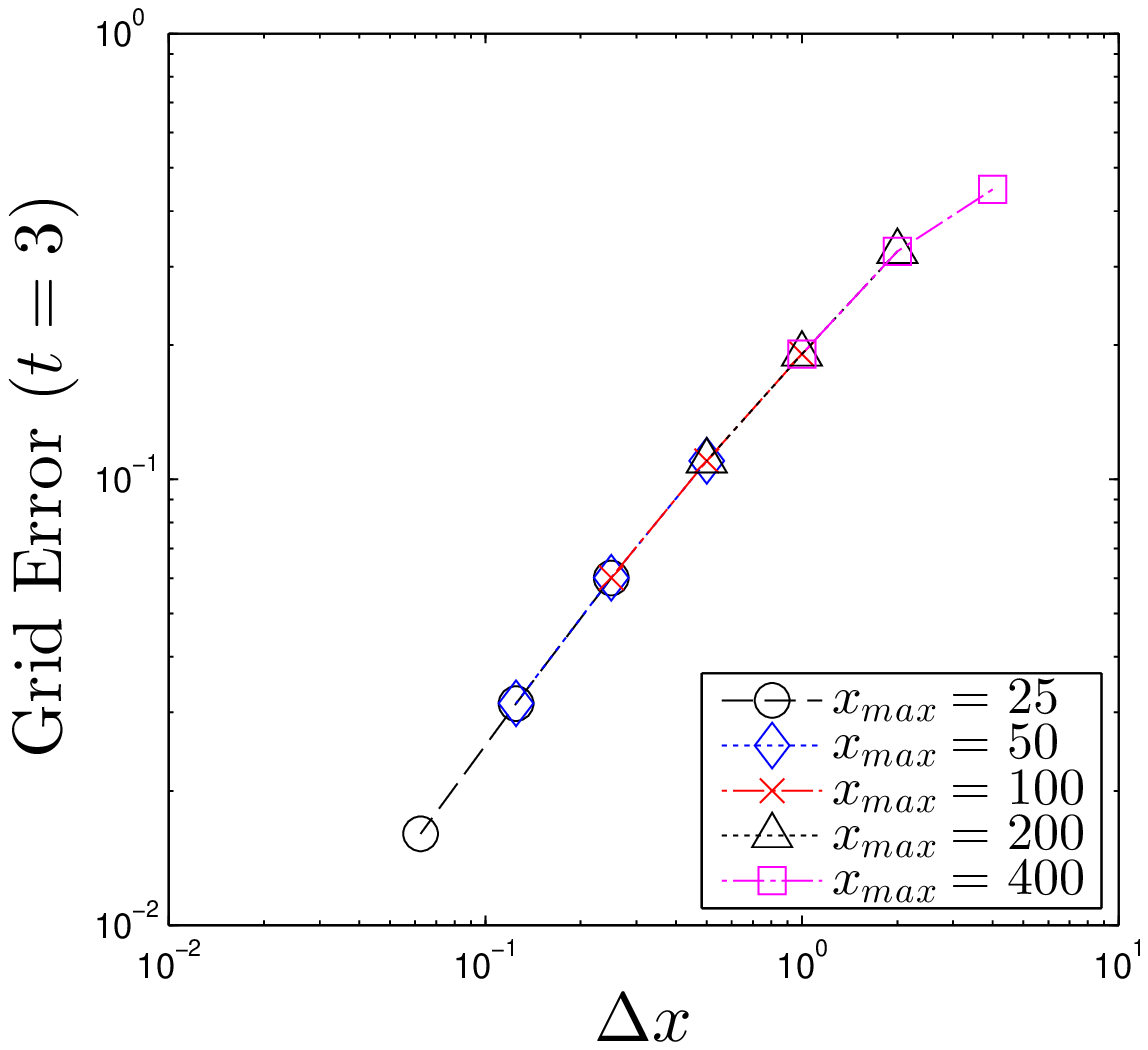}}\hfill{}\subfloat[FLFM]{\centering{}\includegraphics[scale=0.45]{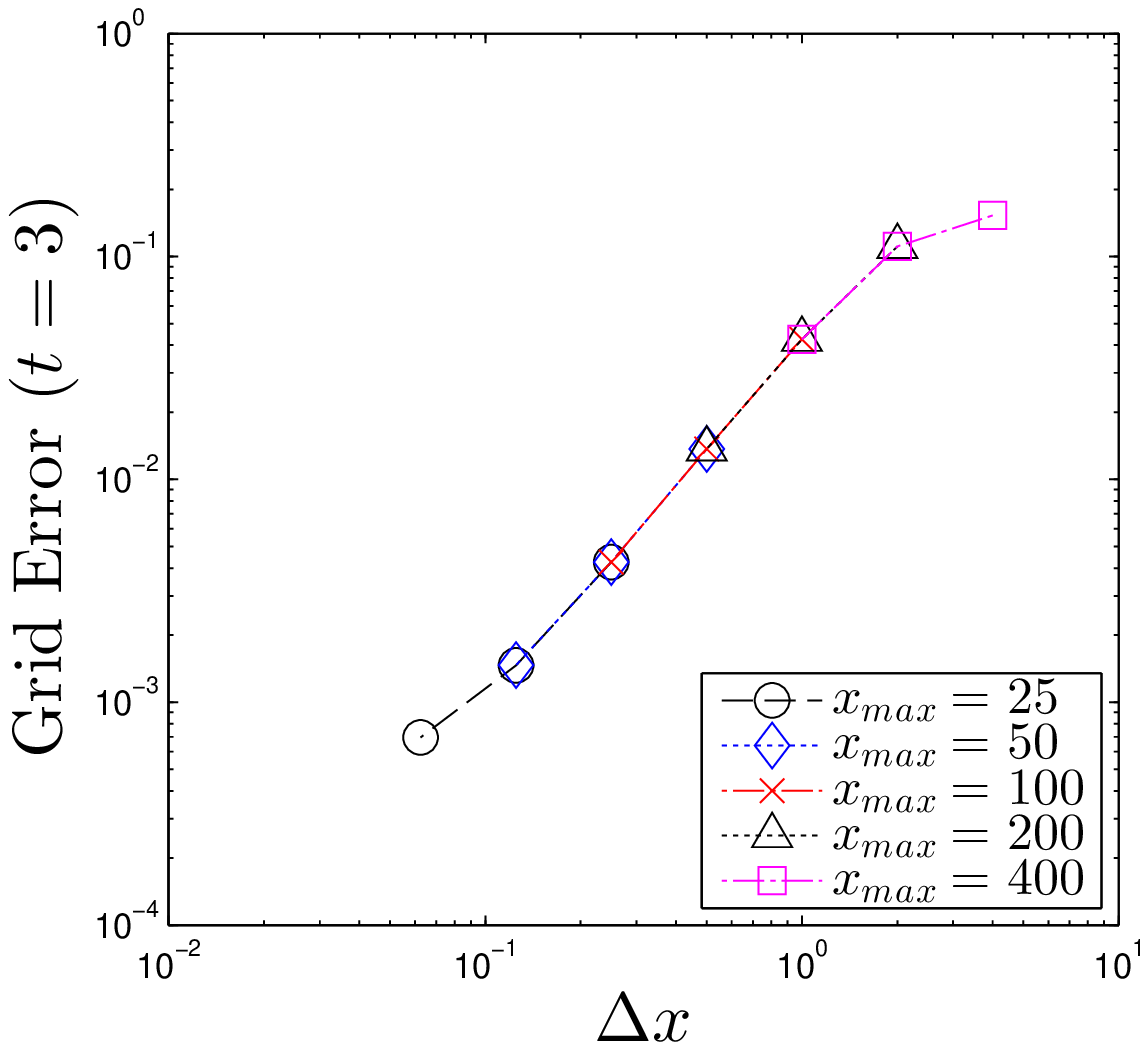}}\caption{Error vs~$\Delta x$ - Grid error at $t=3$ vs.~uniform grid spacing,
$\Delta x$, when $K_{A}(x,y)\equiv1$, for $x_{max}=25,$ $x_{max}=50,$
$x_{max}=100$, $x_{max}=200$, and $x_{max}=400$ in log scale. Clearly,
with a slowly aggregating system, grid spacing plays the primary role
in accuracy for both methods studied in the paper.}
\label{Rdelxplotsais1}
\end{figure}

\begin{figure}[H]
\subfloat[FEM]{\label{error vs delx FEM aisxy}

\centering{}\includegraphics[scale=0.45]{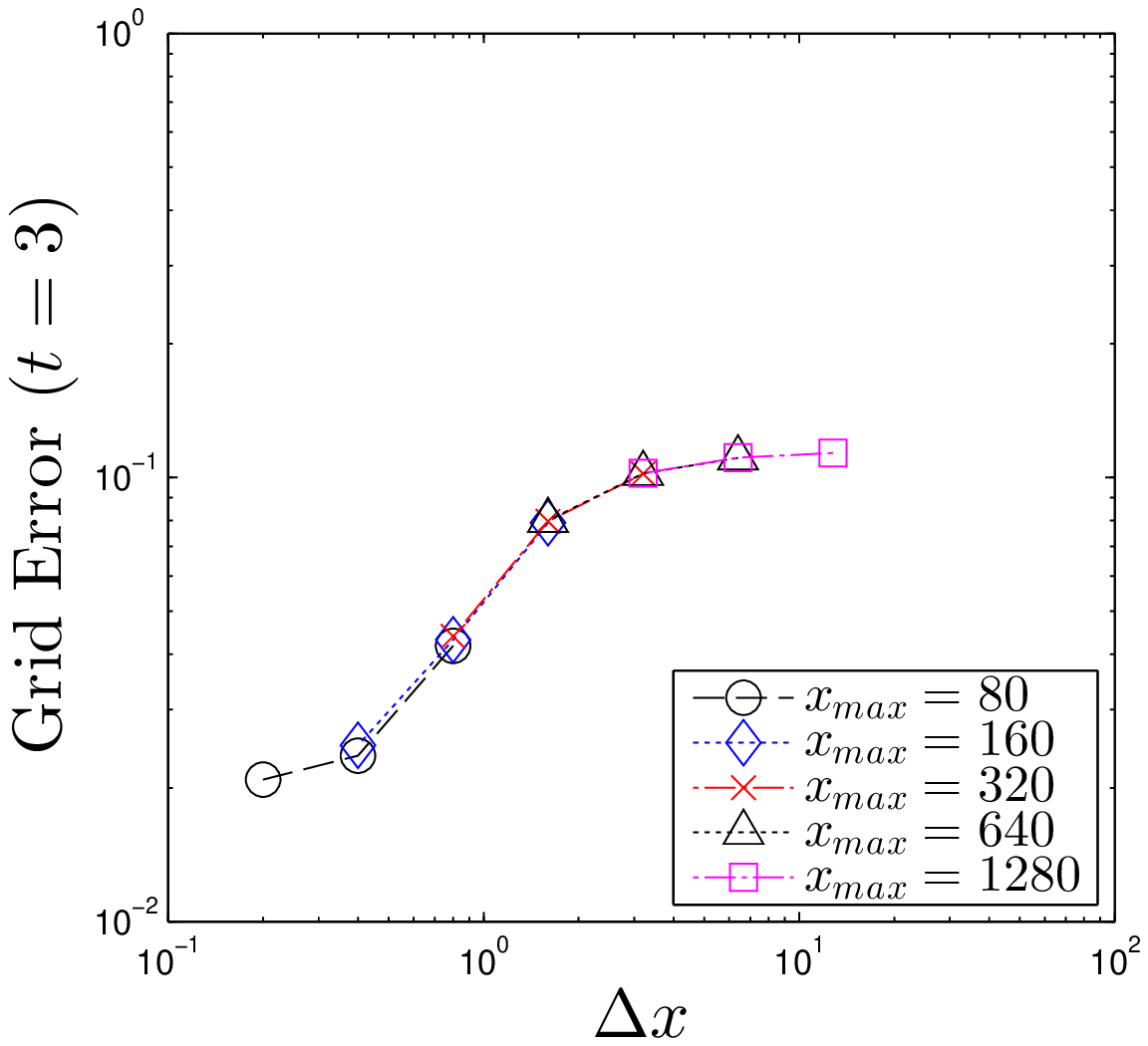}}\hfill{}\subfloat[FLFM]{\label{error vs delx flux aisxy}

\centering{}\includegraphics[scale=0.45]{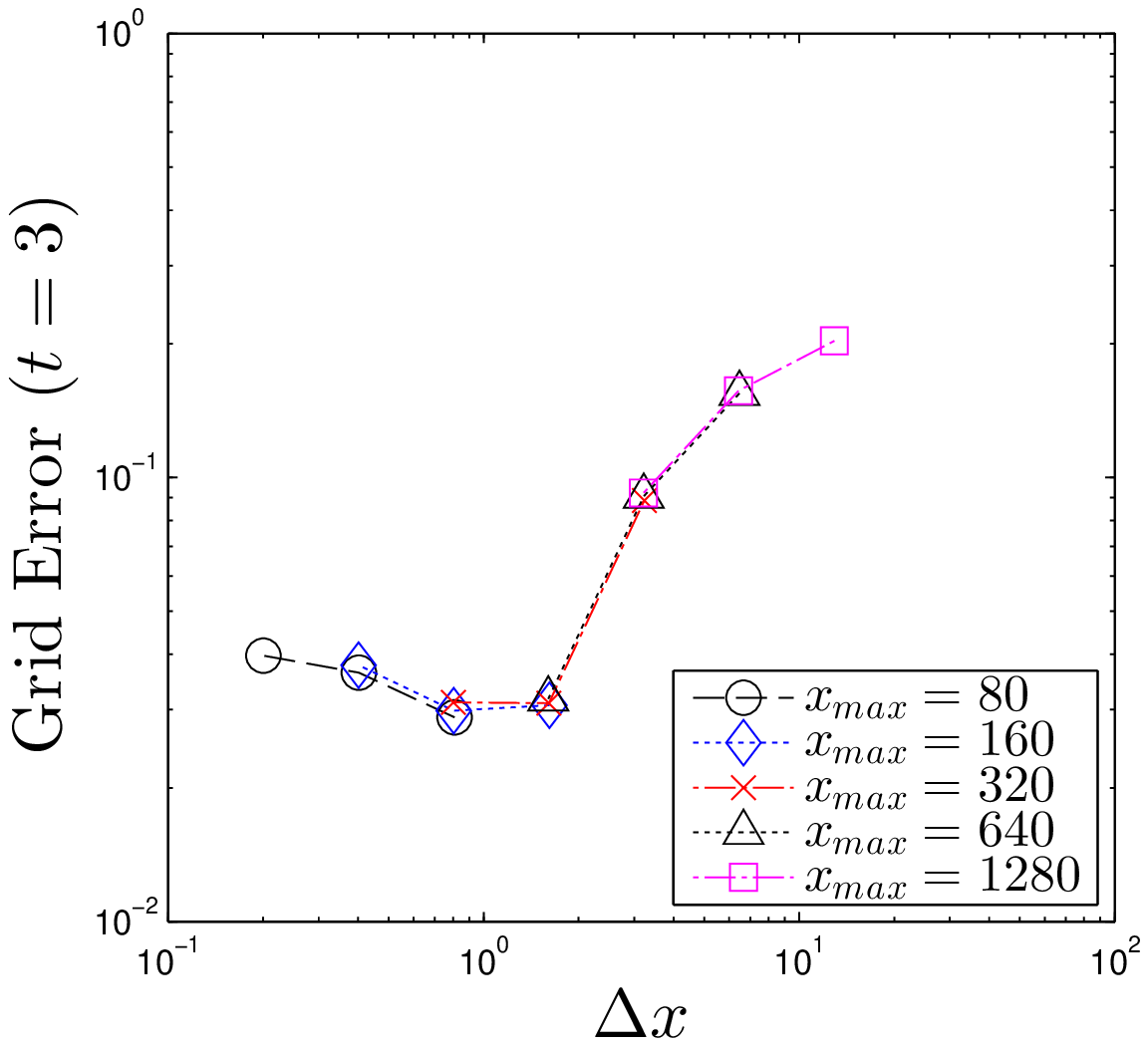}}\caption{Error vs~$\Delta x$ - Grid error at $t=3$ vs.~uniform grid spacing,
$\Delta x$, when $K_{A}(x,y)=xy$, for $x_{max}=80,$ $x_{max}=160,$
$x_{max}=320$, $x_{max}=640$, and $x_{max}=1280$ in log scale.
For the FEM, grid spacing plays a more important role in accuracy
than the truncation parameter, $x_{max}$, but the FLFM experiences
poor behavior when $x_{max}$ is too small.}
\label{errvsdelx aisxy}
\end{figure}

Note that for each run and a given value of $x_{max}$, we initialize
our grid with $x_{1}=0.001$ when $K_{A}(x,y)\equiv1$ and with $x_{1}=0.75$
when $K_{A}(x,y)=xy$. Notice in Figures \ref{Rdelxplotsais1} and
\ref{errvsdelx aisxy}, for a given $\Delta x$, using both the FEM
and the FLFM, the error is the same regardless of the value of $x_{max}$
for both kernels, and we achieve the expected behavior of decreasing
error with refinement of the grid. We suggest that when $K_{A}(x,y)\equiv1$
with the initial conditions used in this study, the analytic solutions
decay so quickly that the numeric approximations do not suffer greatly
from the truncation parameter, $x_{max}$, for either method as evidenced
in Figure \ref{Rdelxplotsais1}. However, when $K_{A}(x,y)=xy$, we
cannot make the same general statement. Using the FLFM with $x_{max}=80$,
$x_{max}=160$, and $x_{max}=320$, we actually achieve reduced error
for larger grid spacing as depicted in Figure \ref{error vs delx flux aisxy},
which is consistent with the poor results noted by Filbet and Laurençot
in \citep{Filbet2004}. Practically speaking, if we wanted to model
a system experiencing rapid aggregation and faced limits on gathering
data for aggregates with large volume, the FEM would provide more
accurate results.

\section{Conclusions and Future Work\label{sec:Conclusions-and-Future}}

The Smoluchowski coagulation equation provides a useful model of particles
in suspension in diverse fields of study. Because only a few analytical
solutions exist, researchers have developed numerical approaches to
approximate solutions to the model. With an eventual goal of comparing
simulated solutions with experimental data, we compare two method's
accuracies in approximating known and fine grid solutions as well
as their accuracies in approximating the zeroth and first moments.
We also compare the computation cost of the two methods.

In \citep{AcklehFitzpatrick1997}, Ackleh and Fitzpatrick report first
order convergence of the FEM in $L^{\infty}[\boldsymbol{X}]$, and
in \citep{Filbet2004}, Filbet and Laurençot, report second order
convergence of the FLFM in $L^{1}[\boldsymbol{X}]$. Our results support
(our conjecture) of second order convergence of the FEM in $L^{1}[\boldsymbol{X}]$
and second order convergence of the FLFM in $L^{1}[\boldsymbol{X}]$,
which eliminates any speed advantage the FLFM was previously understood
to possess. Furthermore, when approximating a fine grid solution the
two methods achieve nearly equal accuracy, with the FLFM achieving
slightly higher accuracy for the multiplicative kernel.

Experimental data comes in different forms such as partial zeroth
moment distributions or partial first moment distributions. We also
theoretically consider the moment approximations and provide numerical
evidence that .the FLFM is slightly more accurate approximating the
zeroth moment for quickly aggregating systems, and the FEM is much
more accurate approximating the zeroth moment when the system aggregates
slowly. In terms of approximating the first moment, the FLFM is more
accurate for both quickly and slowly aggregating systems.

We also identified discretization resolutions where for the same level
of accuracy, the FEM exhibits substantial savings in computational
cost over the FLFM For example, the FEM on 800 grid points offers
an opportunity of computation cost savings over using the FLFM on
100 grid points while achieving nearly equal accuracy. 

We also study the error as a function of both $x_{max}$ and $\Delta x$
in an effort to determine the extent to which the truncation parameter,
$x_{max}$, impacts the overall accuracy of the two schemes. Our initial
results suggest two important results. The first is not surprising:
in general, for the aggregation kernels in this study, reduced grid
spacing plays the predominant role in improving numerical accuracy.
Second, the FLFM suffers from sensitivity to truncation parameter,
$x_{max},$ when $K_{A}(x,y)=xy,$ so if experimental data includes
a small $x_{max}$, the FEM should be the method of choice. These
intriguing patterns will motivate a more extensive theoretical study
of these issue in the future.

\section{Acknowledgements}

This work was supported in part by the National Science Foundation
grant DMS-1225878.

\bibliographystyle{abbrv}
\bibliography{mathbioCU}

\appendix

\section{Example Calculation of Flux Derivative for FLFM\label{sec:flux calcs-1}}

To illustrate how we calculate the right hand side of our second order
spatial approximation in (\ref{eq:flux deriv approx}), we offer the
following example set of calculations. At a given time step, $k$,
consider $J_{3}^{N}(t_{k}).$ We have fixed $r=3$, and then we fix
$x_{mid(1)}$, therefore on a uniform grid, $\tilde{x}\ge x_{3}-x_{mid(1)}=x_{mid(2)}$.
We now know that $x_{mid(1)}+x_{mid(j)}\ge x_{3}\,\mbox{for}\, j\in[2,N]$,
so the contribution to $J_{3}^{N}(t_{k})$ by aggregates in the first
element is{\small 
\[
\Delta xg_{1}^{N}(t_{k})\left[\int_{x_{mid(2)}}^{x_{3}}\frac{K_{A}(x_{mid(1)},y)}{y}g_{2}^{N}(t_{k})dy+\sum_{j=3}^{N-1}\left(\int_{x_{j}}^{x_{j+1}}\frac{K_{A}(x_{mid(1)},y)}{y}g_{j}^{N}(t_{k})dy\right)\right].
\]
}We then consider contributions of aggregates with volume $x_{mid(2)}$,
which leads to $\tilde{x}\ge x_{3}-x_{mid(2)}=x_{mid(1)}$. This implies
$j\in[1,N]$ for contributions to $J_{3}^{N}(t_{k})$ by aggregates
in the second element amounting to{\small 
\[
\Delta xg_{2}^{N}(t_{k})\left[\int_{x_{mid(1)}}^{x_{2}}\frac{K_{A}(x_{mid(1)},y)}{y}g_{1}^{N}(t_{k})dy+\sum_{j=2}^{N-1}\left(\int_{x_{j}}^{x_{j+1}}\frac{K_{A}(x_{mid(1)},y)}{y}g_{j}^{N}(t_{k})dy\right)\right].
\]
}We now have our total flux across $x_{3}$\emph{\footnotesize 
\begin{eqnarray*}
J_{3}^{N}(t_{k}) & = & \Delta xg_{1}^{N}(t_{k})\left[\int_{x_{mid(2)}}^{x_{3}}\frac{K_{A}(x_{mid(1)},y)}{y}g_{2}^{N}(t_{k})dy+\sum_{j=3}^{N-1}\left(\int_{x_{j}}^{x_{j+1}}\frac{K_{A}(x_{mid(1)},y)}{y}g_{j}^{N}(t_{k})dy\right)\right]\\
 & + & \Delta xg_{2}^{N}(t_{k})\left[\int_{x_{mid(1)}}^{x_{2}}\frac{K_{A}(x_{mid(1)},y)}{y}g_{1}^{N}(t_{k})dy+\sum_{j=2}^{N-1}\left(\int_{x_{j}}^{x_{j+1}}\frac{K_{A}(x_{mid(1)},y)}{y}g_{j}^{N}(t_{k})dy\right)\right].
\end{eqnarray*}
}At this point, we can generalize the flux at any given element boundary
($i\in[2,N]$) as{\scriptsize 
\begin{equation}
J_{i}^{N}(t_{k})=\sum_{p=1}^{i-1}\Delta xg_{p}^{N}(t_{k})\left\{ \int_{x_{mid(i-p)}}^{x_{i-p+1}}\frac{K_{A}(x_{mid(p)},y)}{y}dy\, g_{i-p}^{N}(t_{k})+\sum_{j=i-p+1}^{N-1}\int_{x_{j}}^{x_{j+1}}\frac{K_{A}(x_{mid(p)},y)}{y}dy\, g_{j}^{N}(t_{k})\right\} .
\end{equation}
}Specifically, when $K_{A}(x,y)\equiv1$,{\small 
\[
J_{i}^{N}(t_{k})=\sum_{p=1}^{i-1}\Delta xg_{p}^{N}(t_{k})\left\{ \int_{x_{mid(i-p)}}^{x_{i-p+1}}\frac{1}{y}dy\, g_{i-p}^{N}(t_{k})+\sum_{j=i-p+1}^{N-1}\int_{x_{j}}^{x_{j+1}}\frac{1}{y}dy\, g_{j}^{N}(t_{k})\right\} ,
\]
}which after integration gives us{\small 
\[
J_{i}^{N}(t_{k})=\sum_{p=1}^{i-1}\Delta xg_{p}^{N}(t_{k})\left\{ \ln\frac{x_{i-p+1}}{x_{mid(i-p)}}\, g_{i-p}^{N}(t_{k})+\sum_{j=i-p+1}^{N-1}\ln\frac{x_{j+1}}{x_{j}}\, g_{j}^{N}(t_{k})\right\} ,
\]
}and when $K_{A}(x,y)=xy$,{\small 
\[
J_{i}^{N}(t_{k})=\sum_{p=1}^{i-1}\Delta xg_{p}^{N}(t_{k})\left\{ \int_{x_{mid(i-p)}}^{x_{i-p+1}}\frac{x_{mid(p)}y}{y}dy\, g_{i-p}^{N}(t_{k})+\sum_{j=i-p+1}^{N-1}\int_{x_{j}}^{x_{j+1}}\frac{x_{mid(p)}y}{y}dy\, g_{j}^{N}(t_{k})\right\} ,
\]
}which after integration gives us{\small 
\[
J_{i}^{N}(t_{k})=\sum_{p=1}^{i-1}\Delta xg_{p}^{N}(t_{k})\left\{ .5\Delta x\cdot x_{mid(p)}\, g_{i-p}^{N}(t_{k})+\sum_{j=i-p+1}^{N-1}x_{mid(p)}\Delta x\, g_{j}^{N}(t_{k})\right\} .
\]
}{\small \par}

\section{Solution at $t=0$ for $K_{A}(x,y)=xy$\label{sec:t=00003D0 solution notes}}

For the analytical solution used in the paper when $K_{A}(x,y)=xy$,
note that $f(0,x)=\frac{e^{-x}}{x}$, which is not necessarily obvious.
The derivation of those initial conditions follow. From (\ref{eq:xx' solution}),
\begin{equation}
f(t,x)=e^{(-Tx)}\frac{I_{1}(2x\sqrt{t})}{x^{2}\sqrt{t}},
\end{equation}
where
\[
T=\left\{ \begin{array}{ll}
1+t & \mbox{if}\,\, t\le1\\
2\sqrt{t} & \mbox{otherwise}
\end{array}\right.,
\]
and we use the modified Bessel function of the first kind
\[
I_{1}(x)=\frac{1}{\pi}\int_{0}^{\pi}e^{x\,\cos\theta}\cos\theta d\theta.
\]
For this solution, $f(0,x)=\frac{e^{-x}}{x}$. Note that for $t\leq1$,
\[
f(t,x)=\frac{e^{(-x-tx)}I_{1}(2x\sqrt{t})}{x^{2}\sqrt{t}}
\]
with
\[
I_{1}(2x\sqrt{t})=\frac{1}{\pi}\int_{0}^{\pi}e^{(2x\sqrt{t}\cos\theta)}\cos\theta d\theta.
\]
Now note the following:
\begin{enumerate}
\item $\lim_{t\rightarrow0}\frac{I_{1}(2x\sqrt{t})}{x^{2}\sqrt{t}}=\lim_{t\rightarrow0}\frac{1}{\pi}\int_{0}^{\pi}\frac{e^{(2x\sqrt{t}cos\theta)}\cos\theta}{x^{2}\sqrt{t}}d\theta$ 
\item e$^{(2x\sqrt{t}\cos\theta)}=1+(2x\sqrt{t}\cos\theta)+\frac{(2x\sqrt{t}\cos\theta)^{2}}{2!}+\cdots$
\end{enumerate}
therefore,
\begin{eqnarray*}
\lim_{t\rightarrow0}\frac{I_{1}(2x\sqrt{t})}{x^{2}\sqrt{t}} & = & \lim_{t\rightarrow0}\frac{1}{\pi}\int_{0}^{\pi}\frac{\cos\theta+2x\sqrt{t}\cos^{2}\theta}{x^{2}\sqrt{t}}d\theta+\frac{1}{\pi}\int_{0}^{\pi}\lim_{t\rightarrow0}\left[\frac{(2x\sqrt{t}\cos\theta)^{2}\cos\theta}{2!x^{2}\sqrt{t}}+\cdots\right]d\theta\\
 & = & \lim_{t\rightarrow0}\frac{1}{\pi}\int_{0}^{\pi}\frac{\cos\theta}{x^{2}\sqrt{t}}d\theta+\lim_{t\rightarrow0}\frac{1}{\pi}\int_{0}^{\pi}\frac{2\cos^{2}\theta}{x}d\theta+0\\
 & = & 0+\lim_{t\rightarrow0}\frac{2}{\pi x}(\frac{\pi}{2})=\frac{1}{x}
\end{eqnarray*}

from which it follows that
\[
\lim_{t\rightarrow0}f(t,x)=\left[\lim_{t\rightarrow0}e^{(-x-tx)}\right]\left[\lim_{t\rightarrow0}\frac{I_{1}(2x\sqrt{t})}{x^{2}\sqrt{t}}\right]=\frac{e^{-x}}{x}.
\]

\end{document}